\documentclass{amsart}
\usepackage[latin1]{inputenc}
\usepackage[english]{babel}
\usepackage{csquotes}

\usepackage{color} 
\usepackage[colorlinks=true,linktoc=all,linkcolor=blue,citecolor=red,filecolor=pink,urlcolor=blue]{hyperref}
\usepackage[hyperref=true,style=alphabetic,citestyle=alphabetic,backend=bibtex,maxnames=100,doi=false,isbn=false,url=false,giveninits=true]{biblatex}
\bibliography{biblio}

\usepackage{graphicx}
\usepackage{pgf,tikz}
\usetikzlibrary{arrows}
\usetikzlibrary{decorations.markings}

\usepackage{amssymb,amsmath,latexsym}
\theoremstyle{plain}                    
\newtheorem{theorem}{Theorem}[section]
\newtheorem{lemma}[theorem]{Lemma}
\newtheorem{proposition}[theorem]{Proposition}
\newtheorem{corollary}[theorem]{Corollary}
\theoremstyle{definition}
\newtheorem{definition}[theorem]{Definition}
\newtheorem{example}[theorem]{Example}
\newtheorem{question}[theorem]{Question}
\theoremstyle{remark}
\newtheorem{remark}[theorem]{Remark}

\numberwithin{equation}{section}

\newcommand{\cp}{\mathbb{CP}^1}
\newcommand{\pslc}{\textrm{PSL}_2\mathbb{C}}
\newcommand{\slc}{\textrm{SL}_2\mathbb{C}}

\begin{document}

\title[Local deformations of branched projective structures]{Local~deformations of branched~projective~structures: Schiffer~variations and the Teichm\"uller~map}
\author{STEFANO FRANCAVIGLIA}
\address{Dipartimento di Matematica - Universit\`a di Bologna, Piazza di Porta San Donato 5, 40126, 
Bologna, Italy.}
\email{stefano.francaviglia@unibo.it}
\author{LORENZO RUFFONI}
\address{Department of Mathematics - Florida State University, 1017 Academic Way, Tallahassee, FL 
32306-4510, USA.}
\email{lorenzo.ruffoni2@gmail.com}

\keywords{Complex projective structures, movements of branch points, Beltrami differentials, holonomy, hyperelliptic curves.}
\subjclass[2010]{57M50, 14H15, 34Mxx, 32Gxx}
\date{\today}

\begin{abstract}
We study a class of continuous  deformations of branched complex projective structures on closed surfaces of genus 
$g\geq 2$, which preserve the holonomy representation of the structure and the order of the branch points. In the case 
of non-elementary holonomy we show that when the underlying complex structure is infinitesimally preserved the branch 
points are necessarily arranged on a canonical divisor, and we establish a partial converse for hyperelliptic 
structures.
\end{abstract}

\maketitle
\tableofcontents

\section{Introduction}

In this paper we are interested in the deformation theory  of certain classes of complex projective structures on a 
fixed closed surface $S$ of genus $g\geq 2$. These are geometric structures locally modelled on the geometry of the Riemann sphere $\cp$ 
and its group of M\"obius transformations $\pslc$, possibly with branch points; in Thurston's terminology (see \cite{TH97}), they are 
called branched ($\pslc,\cp$)-structures. We refer the reader to Dumas' survey \cite{DU09} for an account of the 
classical theory in the unbranched case, and to Mandelbaum's original papers 
\cite{MA72,MA73} for the branched case. \par

Such a structure induces in particular a complex structure on the underlying surface,  and actually the classical theory 
of Riemann surfaces provides lots of examples. Spherical, Euclidean and hyperbolic geometry have conformal realizations 
inside complex projective geometry, so constant curvature metrics are examples of complex projective structures. Then 
uniformization theory shows that every Riemann surface admits structures of this type, both with and without branch 
points. In particular the deformation space of branched complex projective structures $\mathcal{BP}(S)$ has a natural 
surjective forgetful map to the Teichm\"uller space $\mathcal T(S)$ of the surface 
$S$
$$\pi:\mathcal{BP}(S)\to \mathcal T(S)$$

As usual with geometric structures (see \cite{TH97,DU09}), analytic continuation of local charts gives rise to a global 
invariant known as the holonomy of the structure, which in our context takes the form of a conjugacy class of 
representations $\rho:\pi_1(S)\to \pslc$. We obtain therefore another forgetful map, known as the holonomy map, with 
values in the $\pslc$-character variety $\chi (S)$ of the fundamental group of $S$
$$hol:\mathcal{BP}(S)\to \chi (S)$$

Since the pioneering work of  Poincar\'e on second-order linear ODEs on complex domains (see for instance \cite{PO84}), 
it has been observed that the ratio of two independent solutions of such an equation can be used to define local charts 
for a complex projective structure, which has linked the study of these structures to the classical theory of ODEs, 
especially to the study of questions about their monodromy groups.
In more recent days this has been popularized by the work of  Gallo-Kapovich-Marden in \cite{GKM00} on the Schwarz 
equation on closed surfaces of genus $g\geq 2$, which is the equation locally given by
$$ \left( \dfrac{u''}{u'} \right)' - \dfrac{1}{2} \left( \dfrac{u''}{u'} \right)^2 = q$$
where $q$ is a holomorphic quadratic differential.
By constructing projective structures with specified holonomies,  they show that a non-elementary 
representation $\rho:\pi_1(S)\to \pslc$ arises as the monodromy group of such an equation  if and only if it lifts to $\slc$. In particular it follows from their result that the above map $hol$ is 
essentially surjective, and that a lot of non-discrete representations arise as holonomies of complex projective structures. \par
Our main motivation for this work lies in a more recent paper by Calsamiglia-Deroin-Heu-Loray (see \cite{CDHL16}) in 
which they study $\mathfrak{sl}_2$-systems on  a closed surface of genus $g\geq 2$, i.e. systems of first order ODEs of the form
$$\left( \begin{array}{c} u'_1 \\ u'_2 \\ \end{array} \right)  = 
\left( \begin{array}{cc} a_{11} & a_{12} \\ a_{21} & -a_{11} \\ \end{array} \right) 
\left( \begin{array}{c} u_1 \\ u_2 \\ \end{array} \right) $$
where $a_{ij}$ are holomorphic 1-forms.  In their work they have successfully used the theory of complex projective 
structures to obtain new results about the monodromy groups of these equations, establishing a correspondence between 
$\mathfrak{sl}_2$-systems and certain deformation spaces of complex projective structures.\par

These applications turn the spotlight on the deformation theory of branched complex projective structures  with 
a specified holonomy representation (and possibly a specified number of branch points), and these are the spaces on 
which we focus our attention in this work. Let us denote by $\mathcal M_{k,\rho}$ the deformation space of structures 
having holonomy $\rho:\pi_1(S)\to \pslc$ and $k$ branch points counted with multiplicity.\par

The classical (i.e. unbranched) case corresponds to $k=0$. These spaces $\mathcal M_{0,\rho}$ are known  to be discrete 
by a classical result of Hejhal (see \cite{HE75}). More explicitly, when $\rho$ is a Fuchsian (i.e. discrete, injective 
and real) representation Goldman has shown in \cite{GO87} that every structure in $\mathcal M_{0,\rho}$ is obtained from 
the hyperbolic structure $\mathbb H^2 / \rho(\pi_1(S))$ by a $2\pi$-grafting along a multicurve, a 
surgery introduced originally by Maskit in \cite{MAS69}, and later studied in \cite{GO87, CDF14b}. An analogous result 
has recently been obtained by Baba for generic 
holonomy representations (see \cite{BA12,BA15,BA17}).\par

The branched case (i.e. $k>0)$ requires the introduction of an additional surgery, known as bubbling and  introduced in 
\cite{GKM00}. In the case of Fuchsian holonomy, some results have been obtained which are analogous to Goldman's 
result (see \cite{CDF14, RU19, RU19b}). Moreover the presence of branch points allows certain local deformations 
at branch points, which should be thought as analogous to the classical Schiffer variations for Riemann surfaces (see 
\cite{N}), and have 
actually been shown to provide holomorphic coordinates on these deformation spaces: indeed they have been used in 
\cite[Appendix]{CDF14} to show that if  $\rho$ is a non-elementary representation then $\mathcal M_{k,\rho}$ carries a 
canonical 
structure of $k$-dimensional complex manifold, locally modelled on Hurwitz spaces, and for which the projection to 
Teichm\"uller space $\pi:\mathcal M_{k,\rho}\to \mathcal T(S)$ is holomorphic.\par

A natural question is to better  understand the complex geometry of these deformation spaces, which is of course best 
studied in terms of this projection. This is moreover motivated by the aforementioned correspondence proved in 
\cite{CDHL16} between $\mathfrak{ sl}_2$-systems with irreducible monodromy $\rho$ on a Riemann surface of genus 
$g\geq 2$ and rational curves in $\mathcal M_{2g-2,\rho}$; more precisely they showed that if  $\pi: \mathcal 
M_{2g-2,\rho}\to \mathcal T(S)$ has a fiber with at least three points, then that fiber is actually a rational curve and 
$\rho$ is the monodromy of some $\mathfrak{ sl}_2$-system.\par
This should be contrasted with the unbranched case, in 
which the projection $\pi: \mathcal M_{0,\rho}\to \mathcal T(S)$ is known to be injective by a classical result of 
Poincar\'e (see \cite{PO84}). Recently Baba and Gupta have shown that further projecting to moduli space 
(i.e. forgetting the marking) results in a dense subset of the moduli space of Riemann surfaces (see \cite{BG15}).\par
 
In the branched case an analogous injectivity result has been obtained by Hejhal in \cite[Theorem 15]{HE75} for the 
sub-canonical range $k<2g-2$. The goal of this work is to improve Hejhal's result by studying the infinitesimal 
behaviour of the forgetful map to Teichm\"uller space. Our main result is the following (see Theorem 
\ref{thm_holomoves_can} below); here if $\lambda$ is a partition of $k$ we denote by $\mathcal{M}_{\lambda,\rho}$ the 
deformation space of complex projective structures with holonomy $\rho$ and $k$ branch points counted with multiplicity 
and arranged according to the partition $\lambda$, and we denote by $\pi^\lambda$ the restriction to the stratum
$\mathcal{M}_{\lambda,\rho}$  of the forgetful map $\pi:\mathcal{M}_{k,\rho}\to \mathcal T(S)$.
\begin{theorem}\label{thm_intro1}
Let $g\geq 2$, $\rho:\pi_1(S)\to\pslc$ be non-elementary, $k\leq 2g-2$, $\lambda$ be a partition of 
$k$ and $\sigma \in \mathcal{M}_{\lambda,\rho}$.
If $\sigma$ is a critical point for $\pi^{\lambda}$, 
then the branch points of $\sigma$ form a canonical divisor on the underlying Riemann surface.
\end{theorem}
In particular we get that $k$ must be equal to the degree of a canonical divisor, i.e. $k=2g-2$.
It can be observed by a direct computation that for all the structures on a rational curve in $\mathcal 
M_{2g-2,\rho}$ defined by an $\mathfrak{sl}_2$-system (through the correspondence from \cite{CDHL16}) the branch points are arranged on a canonical divisor on the 
underlying Riemann surface. So the above theorem can be thought as a generalization of this phenomenon.
Motivated once again by the work in \cite{CDHL16}, one would like to find conditions under which some sort of converse 
holds. One could ask for instance if having a canonical branching divisor is actually equivalent to being part of a 
positive-dimensional fiber of the projection to Teichm\"uller space. In order to approach this question we propose a 
study of hyperelliptic structures, i.e. structures endowed with a projective hyperelliptic involution.
Our main result in this direction is the next one (see Theorem \ref{thm_hyp_converse} below for a detailed statement), which is a partial 
converse to Theorem \ref{thm_intro1} in the hyperelliptic setting.
\begin{theorem}\label{thm_intro2}
Let $g\geq 2$, $\rho:\pi_1(S)\to\pslc$ be non-elementary, and let $\sigma \in  \mathcal{M}_{(1,\dots,1),\rho}\subseteq \mathcal{M}_{2g-2,\rho}$
be hyperelliptic. If the branch points of $\sigma$ form a canonical divisor on the underlying Riemann surface, then  
$\sigma$ is a critical point for $\pi$.
\end{theorem}
The proof actually computes local equations for the critical direction, which shows that it is 1-dimensional, and transverse to a natural $(g-1)$-dimensional family of deformations of $\sigma$ through hyperelliptic structures. As a result, for every non-elementary representation which admits hyperelliptic structures (e.g. any non-elementary representation in genus two), we obtain examples of a 1-dimensional family of structures such that all the structures on it have branch points on a canonical divisor and are 
critical points for the projection to Teichm\"uller space, but the family itself is not   a fiber of the 
projection (as already remarked above, it has been shown in \cite{CDHL16} that fibers of the projection have at most dimension 1, in which case they are rational curves associated to $\mathfrak{sl}_2$-systems). \par

Both theorems rely on an explicit computation of the Beltrami differential induced by a movement of branch points; 
these are the aforementioned local deformations at branch points that should be thought as a sort of Schiffer variation 
in this context, as they provide local holomorphic coordinates for $\mathcal{M}_{k,\rho}$ (see 
\cite[Appendix]{CDF14}).\par

This paper is organized as follows. In Section \ref{sec_preliminaries} we introduce branched complex projective 
structures and the deformation spaces we are interested in; moreover in Section \ref{sec_fromODEtoBPS} we review the 
main motivating example from  \cite{CDHL16} for our investigation, namely the correspondence between $\mathfrak{ 
sl}_2$-systems and rational curves in the deformation spaces of projective structures. We include an elementary proof 
that structures coming from this construction have branch points on a canonical divisor. In Section
\ref{sec_complexanalytic} we formally introduce movements of branch points, compute a formula for the Beltrami 
differential induced by such a deformation, and a formula for the contraction of its derivative with a holomorphic 
quadratic differential. This section contains the proof of Theorem \ref{thm_intro1}. Finally 
in Section \ref{sec_hyperelliptic} we introduce hyperelliptic projective structures, construct some explicit examples 
thereof, and provide a proof of Theorem \ref{thm_intro2}.

\vspace{.5cm}
\textbf{Acknowledgements}: We would like to thank Bertrand Deroin and Gabriele Mondello for many useful conversations. 
This work has been partially supported by the European Union's Horizon 2020 research and innovation program under the
Marie Sk\l{}odowska-Curie grant agreement No 777822 (``Geometric and Harmonic Analysis with Interdisciplinary Applications'').

\section{Preliminaries}\label{sec_preliminaries}
Throughout this paper $S$ will be a closed, connected and oriented surface of genus $g\geq 2$, and $\pi_1(S)$ its fundamental group.
We will consider the Riemann sphere $\cp=\mathbb C \cup \{\infty\}$ and its group of holomorphic automorphisms $\pslc$ acting by M\"obius transformations.
We will denote by  $\mathcal{T}(S)$ the Teichm\"uller space of $S$, i.e. the space of marked complex structures on $S$ up to isotopies, and by 
$\chi(S)$ the $\pslc$-character variety of $\pi_1(S)$, i.e. the space of representations into $\pslc$ up to conjugation  by $\pslc$
$$\chi(S)=Hom(\pi_1(S),\pslc)//\pslc$$
where the quotient is taken in the sense of GIT.
By classical results (see \cite{DU09} and references therein) it is known that $\mathcal T(S)$ admits the structure of a smooth connected complex manifold of 
dimension $3g-3$, whereas $\chi(S)$ admits the structure of a complex affine algebraic variety of dimension $6g-6$. In 
this paper  we are interested in the study of geometric structures which have a well-defined underlying complex 
structure and induce $\pslc$-representations of the fundamental group; as such, their deformation spaces come equipped 
with natural maps to $\mathcal T(S)$ an $\chi(S)$, and we will be interested in understanding their relative behaviour.

\subsection{Branched complex projective structures}\label{sec_BPS}
Complex projective structures on $S$ are geometric structures locally modelled on the geometry of the Riemann sphere $\cp$ preserved by the group of M\"obius transformations $\pslc$. In the language of Ehresmann-Thurston geometric structures they are $(\pslc, \cp)$-struc\-tures; we refer the reader to Dumas' survey \cite{DU09} for the classical theory. Branched complex projective structure were introduced by Mandelbaum in \cite{MA72} as a generalization of classic complex projective structures in which branch points are allowed. Here we review the main definitions in this more general context for the convenience of the reader, and to set up notation and terminology.

\begin{definition}\label{def_BPS_charts}
A branched complex projective chart on $S$ is a a pair $(A, d)$ of an open set $A\subseteq S$ and an oriented finite branched cover $d:A\to d(A)\subseteq \cp$ having isolated branch points. Two charts $(A,d), (A',d')$ such that $A\cap A'\neq \varnothing $ are compatible if there exists some $g \in \pslc$ such that $d'=g\circ d$ on $A\cap A'$. A branched complex projective atlas is a covering of $S$ by compatible branched complex projective charts.
\end{definition}
In the following we will call a branched complex projective chart (respectively atlas) just a projective chart  (respectively atlas) for simplicity, whenever no confusion occurs.
\begin{definition}\label{def_BPS_stuff}
A branched complex projective structure (BPS in the following) on $S$ is the datum of a maximal branched complex projective atlas.
A diffeomorphism between two surfaces endowed with BPSs is said to be projective if it is given by M\"obius transformations in local projective charts. 
We denote by $\mathcal{BP}(S)$ the deformation space of BPSs on $S$, namely the space of all marked BPSs on $S$, where two structures are considered to be equivalent if there is a projective diffeomorphism isotopic to the identity between them.
\end{definition}
We refer to Section \ref{sec_holo_fibers} below for an easy way to introduce a natural topology on this space in terms of developing-holonomy pairs.
\begin{example}
Hyperbolic structures are a major source of classical (i.e. unbranched) BPSs, but more general examples of BPSs  include spherical, Euclidean, or hyperbolic metrics of constant curvature with cone points of angle $2\pi n, n \in \mathbb N$, as well as the ``trivial'' BPS defined by any non-constant meromorphic function on a Riemann surface.
\end{example}

\begin{definition}\label{def_order}
A branch point for a BPS $\sigma$ on $S$ is a point $p \in S$ at which local charts have a critical point. 
If a local chart at $p$ has order $n+1$ for some $n\in \mathbb N$, then we say that $p$ has order $ord_\sigma(p)=n$. 
The branching divisor of $\sigma$ is defined to be the divisor $div(\sigma)=\underset{p\in S}{\sum} ord_\sigma (p) p$, and its degree is $|div(\sigma)|=\underset{p\in S}{\sum} ord_\sigma (p)$.
\end{definition}
Notice that the order of a branch point is actually independent of the choice of the chart and that branch points are 
isolated, hence, by compactness of $S$, there is only a finite number of them.

\begin{remark}\label{rmk_teich}
Given a branched cover $\varphi:A\to U$ of an open set $U\subseteq \cp$ with isolated branch points, there is a unique complex structure on $A$ that makes $\varphi$ holomorphic, thanks to Riemann extension theorem. Moreover the changes of coordinates in an atlas for a BPS are given by restrictions of M\"obius transformations, which are in particular holomorphic. As a result, we get that  compatible projective charts induce in particular compatible complex charts, i.e. any BPS has an underlying complex structure, hence we get a natural forgetful map to Teichm\"uller space
$$\pi:\mathcal BP(S)\to \mathcal T(S)$$
which we call the Teichm\"uller map. The fiber of this map over a given $X \in \mathcal T(S)$ can be described in terms of meromorphic quadratic differentials on $X$ with at worst double poles (plus some integrability condition at poles). In this description branch points are encoded in double poles,  and unbranched structures correspond to holomorphic quadratic differentials; more precisely a branch point of order $n-1$ corresponds to a double pole with residue $\frac{1-n^2}{2}$. We refer to \cite{DU09,MA72} for more details.
\end{remark}

\subsection{Holonomy fibers}\label{sec_holo_fibers}
As usual with geometric structures in the sense of Ehresmann and Thurston (see  \cite{DU09,TH97} for more details), local charts can be analytically continued along paths to obtain global objects encoding the structure; for BPSs these are
\begin{itemize}
\item a  representation $\rho:\pi_1(S)\to \pslc$ (called a holonomy representation)
\item  a $\rho$-equivariant orientation preserving smooth map $dev:\widetilde S \to \cp$, which is  which is immersive outside a discrete set of isolated branch points (called a developing map).
\end{itemize}
The equivariance property means that
$$dev\circ \gamma = \rho(\gamma) \circ dev \textrm{ \ for all \ } \gamma \in \pi_1(S)$$
The pair $(dev,\rho)$ is called a developing-holonomy pair, and defines a projective atlas as in the above Definition \ref{def_BPS_charts} just by restriction. 
If $g \in \pslc$ then the pair $(g\circ dev,g\circ \rho \circ g^{-1})$ is considered to be equivalent to $(dev,\rho)$, as they define the same maximal projective atlas; analogously, precomposing $dev$ with a $\pi_1(S)$-equivariant isotopy of $\widetilde S$ results in an equivalent structure. The conjugacy class of holonomy representations is called the holonomy of the structure, and can be recorded by a forgetful map to the $\pslc$-character variety
$$hol:\mathcal{BP}(S)\to \chi(S)$$
which we call the holonomy map.
From this point of view the deformation space $\mathcal{BP}(S)$ can be realized as the space of pairs $(dev,\rho)$, up to the above equivalence relation, and can therefore be equipped with the quotient of the natural compact-open topology on this space of maps.

\begin{remark}
Since a M\"obius transformation $g\in \pslc$ is determined by the image of three points in $\cp$, it is easy to see that a developing map determines uniquely the representation with respect to which it is equivariant; conversely, if $\sigma$ is a BPS whose holonomy has trivial centralizer, then the choice of a holonomy representation $\rho$ in the conjugacy class uniquely determines a unique developing map $dev$ 
such that $(dev,\rho)$ represents the structure $\sigma$. Thanks to these facts, we will use the same notation for a representation and its conjugacy class, and write $\sigma=(dev,\rho)$ to mean that $(dev,\rho)$ is a developing-holonomy pair for $\sigma$, whenever no confusion occurs.
\end{remark}

We are now ready to introduce the deformation spaces we will be mostly concerned with, namely the space of BPSs with prescribed holonomy and number of branch points (counted with multiplicity).
\begin{definition}\label{def_holo_fibers}
Let $\rho \in \chi(S)$. Then we define the holonomy fiber to be
$$\mathcal M_\rho = \{ \sigma \in \mathcal{BP}(S) \ | \ hol(\sigma)=\rho\}$$
We also introduce the following notation: if $k\in \mathbb N$ and $\lambda=(\lambda_1,\dots,\lambda_n)$ is a partition of $k$, then let
$$\mathcal M_{k,\rho} = \{ \sigma \in \mathcal{BP}(S) \ | \ hol(\sigma)=\rho, |div(\sigma)|=k\}$$
$$\mathcal M_{\lambda,\rho} = \{ \sigma \in \mathcal M_{k,\rho} \ | \ \sigma \textrm{ has $n$ branch points, of orders } \lambda_1,\dots,\lambda_n \textrm{ respectively} \}$$
The spaces $\mathcal M_{\lambda,\rho}$ stratify $\mathcal M_{k,\rho}$, and the principal stratum is defined to be the one corresponding to the longest partition, i.e. $\mathcal M_{(1,\dots,1),\rho}$. Structures in the principal stratum are said to be simply branched.
\end{definition}

To investigate non-emptiness of these deformation spaces, we will recall a celebrated theorem by Gallo-Kapovich-Marden (see \cite{GKM00}). Recall that a representation $\rho:\pi_1(S)\to \pslc$ is said to be non-elementary if it has no finite orbits in its action on $\cp \cup \mathbb H^3$. Such a representation 
is in particular irreducible and has trivial centralizer. We also say a representation is liftable if it lifts to a representation into $\slc$. Given our setting, the theorem can be stated as follows.
\begin{theorem}[{\cite[Theorem 1.1.1]{GKM00}}]\label{thm_gkm}
Let $\rho:\pi_1(S)\to \pslc$ be a non-elementary representation. Then
\begin{itemize}
\item $\mathcal M_{0,\rho}\neq \varnothing$   if and only if $\rho$ is liftable.
\item $\mathcal M_{1,\rho}\neq \varnothing$   if and only if $\rho$ is non-liftable.
\end{itemize}
\end{theorem}
Notice that according to our Definition \ref{def_order} a local chart of the form $z\mapsto z^{n+1}$ gives a branch  
point of order $n$, whereas it gives a branch point of order $n+1$ in the notation of \cite[\S 11.2]{GKM00}. 
Moreover it should be noticed that Theorem \ref{thm_gkm} implies that a lot of non-discrete representations arise as 
holonomies of BPSs.\par
Given a structure with a certain holonomy, in order to obtain more BPSs with 
the same holonomy but higher branching, we can perform the following geometric surgery.
\begin{definition}\label{def_bubbling}
Let $\sigma=(dev,\rho)$ be a BPS on $S$. A simple arc $\beta$ is said to be bubbleable for $\sigma$ if the interior of $\beta$ avoids branch points and $dev$ is injective on any lift $\widetilde \beta$ of $\beta$ to the universal cover (we also say that $\beta$ is injectively developed). If $\beta$ is bubbleable on $\sigma$, the bubbling of $\sigma$ along $\beta$ is the structure $Bub(\sigma,\beta)$ obtained by replacing $\beta$ with the disk $B_\beta= \cp \setminus  dev(\widetilde \beta) $, endowed with its natural projective structure (such a disk is called a bubble). 
\end{definition}
\begin{figure}[h]

\begin{center}
\begin{tikzpicture}[scale=0.5]
\node at (-5,0) {$\sigma$};
\draw[xscale=0.8] (-5.5,0) to[out=90,in=180] (-2.75,1.5) to[out=0,in=180] (0,1) to[out=0,in=180] (2.75,1.5) to[out=0,in=90] (5.5,0) ;
\draw[rotate=180,xscale=0.8] (-5.5,0) to[out=90,in=180] (-2.75,1.5) to[out=0,in=180] (0,1) to[out=0,in=180] (2.75,1.5) to[out=0,in=90] (5.5,0) ;

\draw[xshift=-0.3cm,yshift=-0.1cm,xscale=0.8]  (-3.3,0.1) to[out=65,in=180] (-2.5,0.5) to[out=0,in=115] (-1.7,0.1);
\draw[xshift=-0.3cm,yshift=-0.1cm,xscale=0.8]  (-3.4,0.3) to[out=-75,in=180] (-2.5,-0.2) to[out=0,in=-105] (-1.6,0.3);
\draw[xshift=4.3cm,yshift=-0.1cm,xscale=0.8]  (-3.3,0.1) to[out=65,in=180] (-2.5,0.5) to[out=0,in=115] (-1.7,0.1);
\draw[xshift=4.3cm,yshift=-0.1cm,xscale=0.8]  (-3.4,0.3) to[out=-75,in=180] (-2.5,-0.2) to[out=0,in=-105] (-1.6,0.3);

\node at (0.25,0.5) {\color{blue} $\beta$};
\draw[blue,thick] (-1,0) to[out=0,in=180] (-0.5,0.1) to[out=0,in=180] (0.5,-0.1) to[out=0,in=180] (1,0);

\node at (6,0) {$+$};

\node at (13,0) {$\mathbb{CP}^1$};
\node at (10.5,1) {\color{blue} $dev(\widetilde \beta)$};
\draw (8,0) to[out=90,in=180] (10,2) to[out=0,in=90] (12,0) to[out=-90,in=0] (10,-2)  to[out=180,in=-90] (8,0);
\draw (8,0) to[out=-30,in=210] (12,0);
\draw[dashed] (8,0) to[out=30,in=150] (12,0);
\draw[blue,thick] (9,1) to[out=0,in=180] (11,-1);

\node (P) at (4,-2) {};
\node (R) at (4,-4) {};
\path[->,font=\scriptsize,>=angle 90]
(P) edge node[above]{} (R);

\end{tikzpicture}


\begin{tikzpicture}[scale=0.5]
\node at (0,-2.5) {$Bub(\sigma,\beta)$};
\draw[xscale=0.8] (-5.5,0) to[out=90,in=180] (-2.75,1.5) to[out=0,in=180] (0,1.25) to[out=0,in=180] (2.75,1.5) to[out=0,in=90] (5.5,0) ;
\draw[rotate=180,xscale=0.8] (-5.5,0) to[out=90,in=180] (-2.75,1.5) to[out=0,in=180] (0,1.25) to[out=0,in=180] (2.75,1.5) to[out=0,in=90] (5.5,0) ;

\draw[xshift=-0.3cm,yshift=-0.1cm,xscale=0.8]  (-3.3,0.1) to[out=65,in=180] (-2.5,0.5) to[out=0,in=115] (-1.7,0.1);
\draw[xshift=-0.3cm,yshift=-0.1cm,xscale=0.8]  (-3.4,0.3) to[out=-75,in=180] (-2.5,-0.2) to[out=0,in=-105] (-1.6,0.3);
\draw[xshift=4.3cm,yshift=-0.1cm,xscale=0.8]  (-3.3,0.1) to[out=65,in=180] (-2.5,0.5) to[out=0,in=115] (-1.7,0.1);
\draw[xshift=4.3cm,yshift=-0.1cm,xscale=0.8]  (-3.4,0.3) to[out=-75,in=180] (-2.5,-0.2) to[out=0,in=-105] (-1.6,0.3);

\draw[blue,thick] (-1,0) to[out=90,in=180] (0,0.5) to[out=0,in=90] (1,0);
\draw[blue,thick] (-1,0) to[out=-90,in=180] (0,-0.5) to[out=0,in=-90] (1,0);
\node at (-1,0) {$*$};
\node at (1,0) {$*$};

\end{tikzpicture}

\end{center}
 \caption{Bubbling.}
\end{figure}
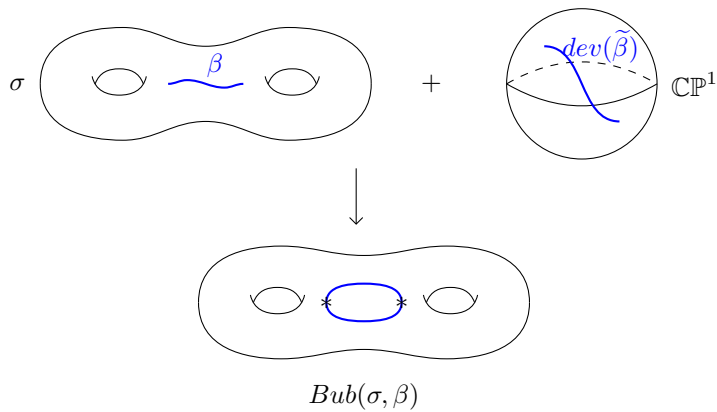
Notice that if $\sigma \in \mathcal{M}_{k,\rho}$ then $Bub(\sigma,\beta)\in \mathcal{M}_{k+2,\rho}$, as two simple branch points appear at the endpoints of $\beta$. This surgery was introduced in \cite{GKM00} as a tool to deform BPSs with the same holonomy, and has later been studied in \cite{CDF14,RU19,RU19b} in the case of quasi-Fuchsian representations. 
These holonomy fibers $\mathcal M_{k,\rho}$ are endowed with natural complex structures as established by the following theorem.
\begin{theorem}[{ \cite[Theorem A.2]{CDF14} }]\label{thm_holofibers_manifolds}
If $\rho$ is non-elementary and $k \in \mathbb N$, then $\mathcal{M}_{k,\rho}$ carries a natural structure of smooth $k$-dimensional complex manifold. Moreover if $\lambda=(\lambda_1,\dots,\lambda_n)$ is any partition of $k$, then $\mathcal{M}_{\lambda,\rho}$ is a smooth $n$-dimensional complex submanifold.
\end{theorem}
The manifold structure is locally modelled on Hurwitz spaces. Naively speaking, local coordinates are obtained by pasting together the local coordinates around the branch points; we refer to the Appendix of \cite{CDF14} and to Section  \ref{sec_surgeries} below for more details about this deformation theory.\par

In this paper we will be mostly concerned with the behavior of these holonomy fibers with respect to the restriction of the Teichm\"uller projection (introduced above in Remark \ref{rmk_teich})
$$\pi:\mathcal M_{k,\rho}\to \mathcal{T}(S)$$
In the classical case of unbranched structures (i.e. when $k=0$), Poincar\'e showed that a complex projective structure is completely determined by its underlying complex structure and its holonomy representation; in other words the projection
$$\pi:\mathcal M_{0,\rho}\to \mathcal{T}(S)$$
is injective. Moreover Baba and Gupta observed in \cite{BG15} the following interesting phenomenon: even though this projection is also proper (by previous work of Tanigawa, see \cite[Theorem 3.2]{TA99}) and its image is therefore a discrete subset of $\mathcal T(S)$, further projecting to moduli space (i.e. forgetting the marking) results in a dense subset (see \cite[Theorem 1.1]{BG15}).\par
In the branched case the dimension of $\mathcal M_{k,\rho}$ is precisely $k$, and therefore the projection can not be injective when $k>3g-3=dim(\mathcal T(S))$ just by dimensional reasons. For small values of $k$ the following injectivity result was proved by Hejhal.
\begin{theorem}[{\cite[Theorem 15]{HE75}}]\label{thm_hejhal}
Let $\sigma \in \mathcal M_{k_1,\rho},\tau \in \mathcal M_{k_2,\rho}$. If $k_1+k_2<4g-4$ and $\pi(\sigma)=\pi(\tau)$, then $\sigma=\tau$
\end{theorem}
As a straightforward consequence the projection
$\pi:\mathcal M_{k,\rho}\to \mathcal{T}(S)$
is injective for $k<2g-2$. The case $k=2g-2$ is of particular interest, as failure of the projection to be injective is deeply related to the existence of solutions of certain ODEs with prescribed monodromy (see next Section \ref{sec_fromODEtoBPS} and \cite{CDHL16} for more details). Our main results concern infinitesimal injectivity of this projection for $k\leq 2g-2$.

\subsection{Rational curves from $\mathfrak{sl}_2$-systems}\label{sec_fromODEtoBPS}
In this  section we review the main motivating example for this paper, i.e. a procedure introduced in  \cite[\S 
6]{CDHL16} that, starting from a certain class of ODEs, constructs a 1-dimensional family of projective structures 
living in the intersection of a Teichm\"uller fiber and a holonomy fiber. We start by introducing the type of ODEs that 
this construction considers; here and in the following we denote by $K_X$ the canonical bundle of a Riemann surface $X$.

\begin{definition}
An $\mathfrak{sl}_2$-system on $S$ is a pair $(X,A)$, where $X\in \mathcal{T}(S)$ and  $A\in H^0(X,K_X\otimes \mathfrak{sl}_2\mathbb C)$.
\end{definition}

Such a datum can naturally be viewed as a linear system of ODEs on a trivial bundle as follows: let $E=X\times \mathbb{C}^2$ the holomorphically trivial rank 2 bundle over $X$; it supports the standard holomorphic connection given by $d$ in local charts, so the choice of $A\in H^0(X,K_X\otimes \mathfrak{sl}_2\mathbb C)$ as above is actually equivalent to the choice of a trace-free holomorphic (hence flat) connection on $E$. Given these data, we can consider the differential equation $d u = A u$ on holomorphic sections of $E$. In local 
holomorphic coordinates such an equation takes the form
$$\left( \begin{array}{c} u'_1 \\ u'_2 \\ \end{array} \right)  = 
\left( \begin{array}{cc} a_{11} & a_{12} \\ a_{21} & -a_{11} \\ \end{array} \right) 
\left( \begin{array}{c} u_1 \\ u_2 \\ \end{array} \right)$$
for some holomorphic coefficients $a_{ij}$. Monodromy of local solutions gives rise to a representation 
$\rho:\pi_1(X)\to \slc$, which we call the monodromy of the system. \par

Notice that both the trivial connection and  the one introduced by $A$ descend to the projectivized bundle $\mathbb 
P(E)=X\times \cp$, hence we get an induced pair of holomorphic connections on it. The flat sections of the former are 
just constant sections $X\times \{c\}$, for $c\in \cp$, whereas the ones of the second give rise to a foliation 
$\mathcal F_A$, generically transverse to the one given by constant sections. Comparing these two foliations defines a 
family of BPSs on $X$ parametrized by $c\in \cp$ as follows: choose a fundamental matrix $\Phi:\widetilde X\to \slc$ for 
the  equation $du=Au$ on the universal cover $\widetilde X $ of $X$, and a point $c\in \cp$; then consider the map
$$\widetilde{X}\to \cp, \quad z \mapsto \Phi(z)^{-1}.c$$
where the action is by M\"obius transformations.\par
This defines a developing map for a BPS on $X$ with holonomy  representation given by the monodromy of the system. We 
denote this structure by $\sigma_{(X,A),c}$.
More concretely, local charts  for it are obtained as follows: fix a point $x_0 \in X $, then flow 
points of $X\times \{c\}$ along leaves of the foliation $\mathcal F_A$ until $\{x_0\} \times \cp$; since the two 
foliations on $\mathbb P (E)$ are generically transverse, this is generically realized by local biholomorphisms from 
open sets in $X\times \{c\}$ to open sets in $\{x_0\}\times \cp$, but branch points arise at tangencies between 
$\mathcal F_A$ and $X\times \{c\}$. It can be computed that precisely $2g-2$ branch points occur, counted with 
multiplicity. We refer to \cite{CDHL16} for more details and a proof of the next statement.
\begin{lemma}
Let $(X,A)$ be an $\mathfrak{sl}_2$-system on $S$ with monodromy $\rho$.
The map $\cp \to \mathcal{M}_{2g-2,\rho},c\mapsto \sigma_{(X,A),c}$ is a holomorphic embedding.

\end{lemma}
We will denote by $\Sigma_{(X,A)}=\{\sigma_{(X,A),c} \ | \ c \in \cp\}  \subset \mathcal{M}_{2g-2,\rho}$ the induced 
family of BPSs.
Notice that all the structures in it have the same underlying complex structure $X$ and the same holonomy $\rho$.
Direct computations then show they also share a certain property about their branching divisor.
\begin{proposition}\label{prop_ODEBPS_can}
 For any $ \sigma_{(X,A),c} \in \Sigma_{(X,A)}$  the branching divisor of $ \sigma_{(X,A),c}$ is a canonical divisor on $X$.
\end{proposition}
\proof  
The branch points of $\sigma_{(X,A),c}$ occur exactly at tangencies between the foliation $\mathcal{F}_A$ defined by 
$A$ and the horizontal section $X\times \{[c_1:c_2]\}$ inside the projectivized bundle $\mathbb P(E)=X\times \cp$, where $c=[c_1:c_2]$.  A 
direct computation using the local expression for the foliations shows that this occurs exactly at points $z\in X$ at 
which $(c_1,c_2)$ is an eigenvector for $A(z)$. Since $E=X\times \mathbb C^2$ is the trivial bundle, we have that 
$H^0(X,K_X\otimes \mathfrak{sl}_2\mathbb C)=H^0(K_X)^{\oplus 3}$; more explicitly, $A$ can be written as
$$A=\left( \begin{array}{cc}
            a_{11} & a_{12} \\
            a_{21} & -a_{11} \\
           \end{array} \right), \textrm{\ for \ } a_{ij} \in H^0(X,K_X)$$
Then we see that $z\in X$ is a point at which  $(c_1,c_2)$ is an eigenvector for $A(z)$ if and only if the following 
conditions are satisfied 
$$\left\lbrace \begin{array}{l}
            c_1a_{11}(z)+c_2a_{12}(z)=c_1\lambda(z) \\
            c_1a_{21}(z)-c_2a_{11}(z)=c_2\lambda(z) \\
           \end{array} \right.$$          
where $\lambda \in H^0(X,K_X)$ is an eigenvalue of $A$. Since $c_1,c_2$ are not both zero, let us assume $c_1\neq 0$, 
express $\lambda(z)=a_{11}(z)+\frac{c_2}{c_1}a_{12}(z)$ and obtain therefore that
$$c_1^2a_{21}(z)-2c_1c_2a_{11}(z)-c_2^2a_{12}(z)=0$$
In other words the branching divisor of $\sigma_{(X,A),c}$ is exactly the zero divisor of the holomorphic 1-form
$\Theta_{A,c}=c_1^2a_{21}-2c_1c_2a_{11}-c_2^2a_{12}$.
\endproof
As mentioned above, by construction the rational curve $\Sigma_{(X,A)}$ lives in the intersection of a Teichm\"uller fiber and a holonomy fiber; in particular all the structures on it are critical points for the projection $\mathcal{ M}_{2g-2,\rho}\to \mathcal T(S)$.
The main result of this paper (see  Theorem \ref{thm_holomoves_can}) is a generalization of the above statement to structures which are critical points for the projection to Teichm\"uller space, but are not necessarily on a fiber thereof.

\section{The local structure of the holonomy fiber}\label{sec_complexanalytic}
In this section we study the local behavior of Teichm\"uller map  along a holonomy fiber, i.e. the natural forgetful map to Teichm\"uller space
$$\pi:\mathcal M_{k,\rho} \to \mathcal{ T}(S)$$
obtained by sending a BPS to its underlying complex structure (see Remark \ref{rmk_teich}).
By Theorem \ref{thm_holofibers_manifolds} above we know that as soon as $\rho$ is non-elementary the space $\mathcal M_{k,\rho}$ carries a natural structure of smooth $k$-dimensional complex manifold, with respect to which $\pi$ is holomorphic. 
Throughout the rest of the paper the holonomy is always assumed to be non-elementary. 
Our aim is to show that a BPS $\sigma$ can be deformed preserving both its holonomy $\rho$ and the underlying complex structure $X$ only if its branching divisor is a canonical divisor on $X$.

\subsection{Movements of branch points}\label{sec_surgeries}

In the previous section (see Definition \ref{def_bubbling}) we have introduced bubbling as a geometric surgery that turns a BPS into another BPS with the same holonomy by increasing the total branching order by 2. In view of Theorem \ref{thm_gkm} this allows in particular to create BPSs with arbitrarily high branching order and prescribed monodromy. In this section we introduce another geometric surgery that preserves both the holonomy and the structure of the branching divisor, therefore providing a way to navigate the spaces $\mathcal{M}_{k,\rho}$.

\begin{definition}\label{def_move}
Let $\sigma$ be a BPS on $S$. Given a branch point $p$ of $\sigma$, a movement of branch point at $p$ is the following 
deformation: pick a local branched projective $(A,d)$ chart at $p$, and replace it with the chart obtained by 
postcomposing $d:A\to \cp$ with a smooth isotopy compactly supported in $d(A)$. If the isotopy moves $d(p)$  to $z\in 
d(A)$ we will denote the resulting structure by $Move(\sigma,z)$ (see Figure \ref{pic_movement}).
\end{definition}
\begin{remark}\label{rmk_move_param}
The result of a movement of branch points depends, up to isomorphism,  only on the final image of 
the branch point, i.e. the point of $d(A)$ to which $d(p)$ moves through the isotopy; in particular 
we can choose the isotopy to be a straight-line isotopy in the chosen  projective chart. Moreover it 
is clear that such a deformation can be performed independently at distinct branch points; if 
$p_1,\dots,p_n$ are the branch points of $\sigma$ and $(A_i,d_i)$ is a local chart at $p_i$, with $A_i\cap A_j=\varnothing$ for $i\neq j$, then 
such a deformation is conveniently parametrized by $\prod_{i=1}^n d_i(A_i)\subseteq \prod_{i=1}^n 
\cp $. In the above notations, we denote by $Move(\sigma,z_1,\dots,z_n)$ the structure obtained by 
moving $d_i(p_i)$ to $z_i$.
\end{remark}
\begin{remark}\label{rmk_move_in_stratum}
Moving branch points does not change the structure of the branching divisor, i.e. the number and order of branch points.
If $\lambda$ is any partition of $k\in\mathbb N$ and  $\sigma \in \mathcal{M}_{\lambda,\rho}$ then $Move(\sigma,z_1,\dots,z_k) \in \mathcal{M}_{\lambda,\rho}$, and it follows from Theorem \ref{thm_holofibers_manifolds} these deformations account for a full neighborhood of $\sigma$ inside $\mathcal{M}_{\lambda,\rho}$. For simply branched structures (i.e. those in the principal stratum, as defined in Definition \ref{def_holo_fibers}), this is also a neighborhood in $\mathcal{M}_{k,\rho}$. For structures in lower strata, a full neighborhood in $\mathcal{M}_{k,\rho}$ should also include deformations which 
split higher order branch points into lower order ones; we refer to \cite[Appendix]{CDF14} for more details.
\end{remark}

\subsection{Beltrami differentials for movements of branch points}
We begin by computing an explicit formula for the Beltrami differential  associated to the 
deformation of the underlying complex structure $X$ induced by a movement of branch points: more 
precisely we will consider a BPS $\sigma$ and a 1-parameter deformation $\{\sigma_t \ | \ t\in 
[0,1]\}$ obtained by moving one branch point in a given direction, and we will compute the 
1-parameter family of Beltrami differentials $\{\mu_t \ | \ t\in [0,1]\}$ induced by it; this is 
done below in Lemma \ref{lem_beltrami}.
We will also compute the first order approximation of this 1-parameter family of Beltrami 
differential, which we will regard as another Beltrami differential approximating the initial 
infinitesimal behavior of $\{\mu_t \ | \ t\in [0,1]\}$.\par

Recall that if $f:X\to Y$ is a smooth map between 
Riemann surfaces, its Beltrami differential is defined to be the smooth section of the line bundle 
$\overline K_X \otimes K_X^{-1}$ given in local complex charts by 
$\frac{\partial f}{\partial \overline{z}} \left( \frac{\partial f}{\partial z} \right)^{-1}$, and is a measure of the failure of $f$ to be holomorphic. We refer to \cite{IT92} for background about Beltrami differentials, and in particular for the existence of a natural duality between Beltrami differentials and holomorphic quadratic differentials, which provide respectively a description of the tangent and cotangent space to Teichm\"uller space.\par

Let $\sigma \in \mathcal{M}_{k,\rho}$, let $X \in \mathcal T(S)$ be the underlying complex structure, and let $p$ be a branch point of order $ord(p)=m-1\geq 1$. Recall from Definition \ref{def_order} that this means that local projective charts at $p$ are branched covers of order $m$. First of all we find a normal local expression for the projective charts.
\begin{lemma}\label{lem_projchartincomplexchart}
 For any projective local chart $(A,d)$ for $\sigma$ at $p$ there exists a complex local chart $(A,\varphi)$ for $X$ at 
$p$ such that $d\circ \varphi^{-1}:\varphi(A)\to d(A)$ is given by $d\circ \varphi^{-1}(z)=z^m+o(z^m)$.
\end{lemma}
\proof Given a complex local chart $(A,\psi)$ and a series expansion for $d\circ \varphi^{-1}$, a translation can be used to fix the constant term, and a dilatation can be used to fix the next non-zero term. In other words $\varphi$ can be obtained by post-composing $\psi$ with a complex affine transformation.\endproof

Let us fix a branched complex projective atlas $\mathcal U$ for $\sigma$. In the notations of Lemma \ref{lem_projchartincomplexchart}, let us pick an open set  $B$ such that $p\in B \subset A$ and such that for any other projective chart $(U,g)\in\mathcal{U}\setminus \{(A,d)\} $ we have $A\cap U \subset A\setminus B$. Let $z$ be a complex coordinate on $\varphi(A)$ and $w$ a projective coordinate on $d(A)$ and let us denote by $c$ the holomorphic map $c=d\circ \varphi^{-1}:\varphi(A)\to d(A),c(z)=z^m+o(z^m)$.\par

\begin{figure}[h]

\begin{center}
\begin{tikzpicture}[scale=0.8]

\fill[gray!20] (-4,0) circle (2cm);
\fill[white] (-4,0) circle (1.5cm);
\fill[gray!20] (-4,0) circle (.5cm);
\draw (-4,0) circle (2cm);
\draw (-4,0) circle (1.5cm);
\draw (-4,0) circle (.5cm);

\begin{scope}[xscale=-1]
\fill[gray!20] (-4,0) circle (2cm);
\fill[white] (-4,0) circle (1.5cm);
\fill[gray!20] (-4.8,0) circle (.5cm);
\draw (-4,0) circle (2cm);
\draw (-4,0) circle (1.5cm);
\draw (-4.8,0) circle (.5cm);
\end{scope}

\begin{scope}[yshift=6.5cm]
\fill[gray!20] (-4,0) circle (2cm);
\fill[white] (-4,0) circle (1.5cm);
\fill[gray!20] (-4,0) circle (.5cm);
\draw (-4,0) circle (2cm);
\draw (-4,0) circle (1.5cm);
\draw (-4,0) circle (.5cm);
\end{scope}

\begin{scope}[xscale=-1,yshift=6.5cm]
\fill[gray!20] (-4,2.5) to[out=0,in=90] (-2.5,1) to[out=-90,in=90] (-3,0) to[out=-90,in=90] 
(-2.5,-1) 
to[out=-90,in=0] (-4,-2.5) to[out=180,in=-90] (-5.5,-1)   to[out=90,in=-90] (-5,0) 
to[out=90,in=-90] 
(-5.5,1) to[out=90,in=180] (-4,2.5);

\draw (-4,2.5) to[out=0,in=90] (-2.5,1) to[out=-90,in=90] (-3,0) to[out=-90,in=90] (-2.5,-1) 
to[out=-90,in=0] (-4,-2.5) to[out=180,in=-90] (-5.5,-1)   to[out=90,in=-90] (-5,0) to[out=90,in=-90] 
(-5.5,1) to[out=90,in=180] (-4,2.5);

\fill[white] (-4,2) to[out=0,in=90] (-3,1) to[out=-90,in=90] (-3.25,0) to[out=-90,in=90] (-3,-1) 
to[out=-90,in=0] (-4,-2) to[out=180,in=-90] (-5,-1)   to[out=90,in=-90] (-4.75,0) 
to[out=90,in=-90] (-5,1) to[out=90,in=180] (-4,2);

\draw (-4,2) to[out=0,in=90] (-3,1) to[out=-90,in=90] (-3.25,0) to[out=-90,in=90] (-3,-1) 
to[out=-90,in=0] (-4,-2) to[out=180,in=-90] (-5,-1)   to[out=90,in=-90] (-4.75,0) 
to[out=90,in=-90] (-5,1) to[out=90,in=180] (-4,2);

\fill[gray!20] (-4,0) circle (.5cm);
 \draw (-4,0) circle (.5cm);
 \end{scope}

\begin{scope}[xshift=4cm,yshift=12cm]
\fill[gray!20] (-4,0) circle (2cm);
\fill[white] (-4,0) circle (1.5cm);
\fill[gray!20] (-4,0) circle (.5cm);
\draw (-4,0) circle (2cm);
\draw (-4,0) circle (1.5cm);
\draw (-4,0) circle (.5cm);
\end{scope}


\draw[->] (-1.5,0) -- (1.5,0);
\draw[->] (-1.5,6) -- (1.5,6);
\draw[->] (-4,4) -- (-4,2.5);
\draw[->] (4,3.5) -- (4,2.5);
\draw[->] (-2,10) -- (-3,9);
\draw[->] (2,10) -- (3,9);
\draw[->] (-2.5,12) -- (-4,12) to[out=180,in=90] (-7,6.5) to[out=-90,in=120] (-6,2) ;
\draw[xscale=-1,->] (-2.5,12) -- (-4,12) to[out=180,in=90] (-7,6.5) to[out=-90,in=120] (-6,2) ;

 \node at (-2.5,13) {$A$};
 \node at (0.5,13) {$B$};
 \node at (-0.5,11.5) {$C$};
 \node at (0,12) {$p$};
 \node at (-6,12) {$d$};
 \node at (-2,9) {$\varphi$};
 \node at (-3,3) {$c=d\varphi^{-1}$};
 \node at (-4,6.5) {$0$};
 \node at (-4,0) {$0$};
 \node at (-6,5) {$z\in \mathbb{C}$};
 \node at (-6.3,-1.5) {$w\in \cp$};
 \node at (0,-0.5) {$H_t(w)=w+tq\eta(w)$};
 \node at (6.4,-1.5) {$v\in \cp$};
  \node at (6,5) {$u\in \mathbb{C}$};
 \node at (3,3) {$c_t=d_t\varphi_t^{-1}$};
 \node at (2,9) {$\varphi_t$};
 \node at (6,12) {$d_t=H_td$};
 \node at (0,6.5) {$F_t$};
 \node at (4,0) {$0$};
 \node at (4,6.5) {$0$};
 \node at (4.8,0) {$tq$};

\end{tikzpicture}
\end{center}

%
%
%
 \caption{Local analysis of the movement of a branch point.}\label{pic_movement}
 \end{figure}
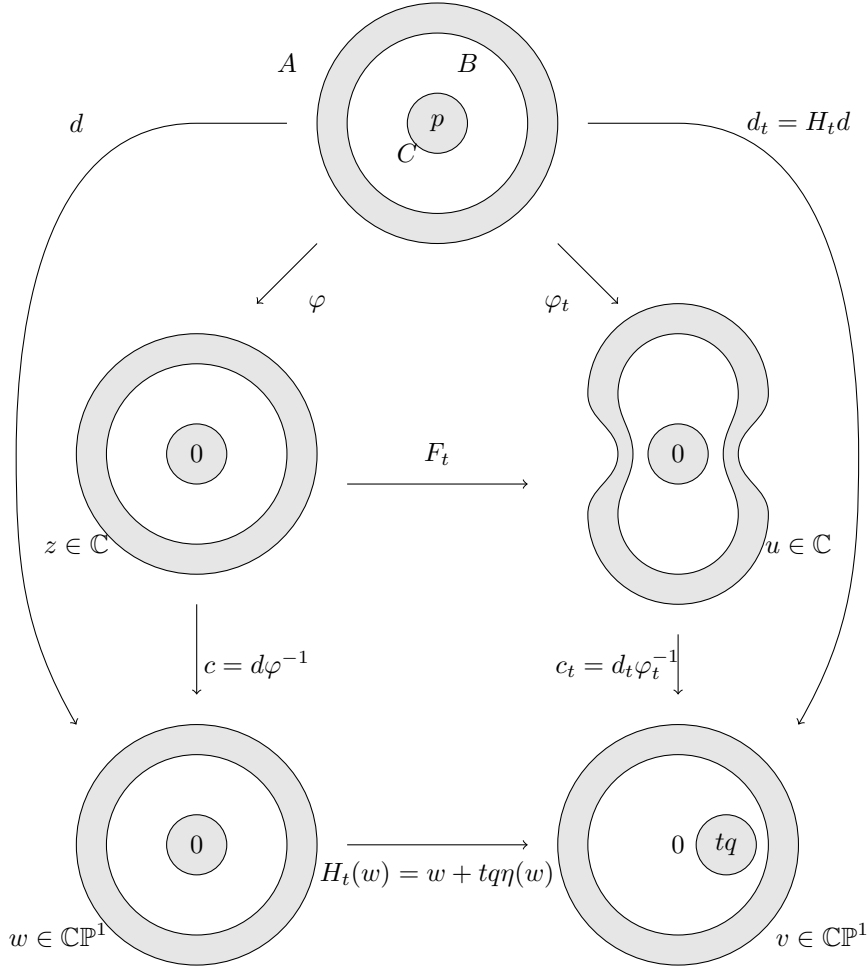

Let us consider a point $q \in d(A)$, an open  neighborhood $C\subset \overline{C}\subset  B$ of $p$, and a smooth bump function $\eta:d(A)\to [0,1]$ such that
\begin{itemize}
\item $\eta$ is compactly supported in $d(B)$
\item $\eta=1$ on $d(C)$
\end{itemize}
Then we get a well-defined isotopy $$H:[0,1]\times d(A)\to d(A)$$ $$ H(t,w)=w+tq\eta(w)$$ supported on $d(B)$. In particular, for any $t\in [0,1]$ the map $H_t:d(A)\to d(A),H_t(w)=H(t,w)$ is a smooth isotopy of $d(A)$ which is projective on $d(A)\setminus d(B)$ and on $d(C)$, where it coincides respectively with the identity and with the translation $w\mapsto w+tq$ (see Figure \ref{pic_movement}). 

Let us define $d_t=H_t\circ d:A\to d(A)$. Then replacing the chart $(A,d)$ with the chart $(A,d_t)$ gives a new BPS on $S$, which we denote by $\sigma_t$. Let us denote by $X_t$ the complex structure underlying $\sigma_t$. With the same argument of Lemma \ref{lem_projchartincomplexchart}, we can obtain a local complex chart $(A,\varphi_t)$ for the complex structure $X_t$ such that $c_t:=d_t\circ \varphi_t^{-1}:\varphi_t(A)\to d_t(A)=d(A)$ is given by $c_t(u)=tq+u^m+o(u^m)$ for a coordinate $u$ on $\varphi_t(A)$. We can then choose a diffeomorphism $F_t:\varphi(A)\to \varphi_t(A)$ which lifts  $H_t:d(A)\to d_t(A)=d(A)$ (i.e. $H_t \circ c = c_t \circ F_t$); among the possible lifts, there is exactly one which satisfies $F_t\circ \varphi = \varphi_t$, and we choose this one.\par
Notice that $F_t$ will be holomorphic on $\varphi(A)\setminus \varphi(B)$ and on $\varphi(C)$, but not elsewhere; also notice that if $\varphi (A) = \varphi_t(A)$ then the deformation of the complex structure would be trivial. Anyway this is never the case: indeed $F_t$ can be regarded as a sort of Schiffer variation of the underlying complex structure (see \cite{N}). To get a feeling of how this works, we can look for a geometric description of the domain $\varphi_t(A)$ as a subset of the complex plane with coordinate $u\in\mathbb{C}$. Without loss of generality, we can assume that $d(A)=\mathbb{D}=\{|w|<1\}$; then $d_t(A)=\{|v|< 1\}$. Then the map $c_t$ will take the form 
$v=c_t(u)=tq+\sum_{k\geq m}e_ku^k$, with $e_m=1$, and the domain $\varphi_t(A)$ of the $u$-plane is the one which is 
mapped by $c_t$ to the unit disk $d_t(A)=\{|v|<1\}$ in the $v$-plane. A direct computation shows that

$\varphi_t(A)$ is a domain bounded by a curve defined by an equation of the form 
$$t^2|q|^2+2t\sum_{k=m}^{2m-1}Re(\overline{q}e_ku^k) +|u|^{2m}+o(u^{2m})-1=0$$

\begin{example}
 Let us consider the base case of a simple branch point (i.e. $m=2$); let us also assume that the local charts take the form $d(A)=\{|w|=1\}$, $d_t(A)=\{|v|=1\}$, $w=c(z)=z^2$ and $v=c_t(u)=tq+u^2$. Then the above computation shows that $\varphi_t(A)$ is the domain of the $u$-plane bounded by the curve
$$ t^2|q|^2+2tRe(\overline{q}u^2)+|u|^4-1=0$$
Plotting this curve for values such that $tq=1-\varepsilon$,  for $\varepsilon>0$ small enough, we obtain  pictures consistent to the one shown on the right in Picture \ref{pic_movement}. 
The map $F_t:\varphi(A)\to \varphi_t(A)$ from the unit disk in the $z$-plane to this bean-like domain is smooth and holomorphic near the boundary and around $0$. We can directly compute that
$$F_t(z)^2+tq=c_t(F_t(z))=H_t(c(z))=H_t(z^2)=z^2+tq\eta(z^2)$$
i.e. the map $F_t:\varphi(A)\to\varphi_t(A)$ looks like $F_t(z)=\sqrt{z^2+tq(\eta(z^2)-1)}$, which coincides with $F_t(z)=z$ around $z=0$ and with $F_t(z)=\sqrt{z^2-tq}$ near the boundary of $\varphi(A)$.
\end{example}

We are now going to compute the Beltrami  differential of the identity map $id_S$ of the surface $S$ 
considered as a map between the marked Riemann surfaces $X$ and $X_t$. For convenience, let us 
introduce the notation $F(t,z)=F_t(z)$ and $\mu(t,z)=\mu_t(z)$.
\begin{lemma}\label{lem_beltrami}
In the above notations, the Beltrami differential $\mu_t$ of $id_S:X\to X_t$ is zero outside $A$, and with respect to the coordinate $z$ over $\varphi(A)$ it is given by the following expression
 $$\mu_t(z)=\mu(t,z)= \dfrac{tq\dfrac{\partial \eta}{\partial 
 \overline{w}}(c(z))\overline{\dfrac{\partial c}{\partial z}(z)}}{\left(1+tq\dfrac{\partial 
 \eta}{\partial w}(c(z))\right)\dfrac{\partial c}{\partial z}(z)} $$
\end{lemma}
\proof
With the atlases we are using, the identity map reads as the identity for any choice of charts, with the exception of the choice of charts $(A,\varphi)$ for $X$ and $(A,\varphi_t)$ for $X_t$; by construction, in these charts it reads as the map $F_t:\varphi(A)\to \varphi_t(A)$, hence we reduce to compute the Beltrami differential of $F_t$. We recall from the above construction the relation $c_t(F(t,z))=H(t,c(z))=c(z)+tq\eta(c(z))$. Taking the derivative with respect to $z$  we obtain
$$\dfrac{\partial c_t}{\partial u}(F(t,z))\dfrac{\partial F}{\partial z}(t,z)=
\dfrac{\partial c}{\partial z}(z)\left( 1+tq\dfrac{\partial \eta}{\partial w}(c(z))\right)$$
and taking the derivative with respect to  $\overline{z}$  we obtain
$$\dfrac{\partial c_t}{\partial u}(F(t,z))\dfrac{\partial F}{\partial \overline{z}}(t,z)=
\overline{\dfrac{\partial c}{\partial z}(z)} tq\dfrac{\partial \eta}{\partial \overline{w}}(c(z)))$$
Comparing the two equalities we get the desired expression for $\mu_t(z)=\dfrac{\dfrac{\partial F}{\partial \overline{z}}(t,z)}{\dfrac{\partial F}{\partial z}(t,z)}$.
\endproof

Let us now compute the first order approximation at $t=0$ for this 
1-parameter family of deformations. 
\begin{lemma}\label{lem_beltrami_firstorder}
 $ \dfrac{\partial \mu}{\partial t}(0,z)=\dfrac{q\dfrac{\partial \eta}{\partial 
 \overline{w}}(c(z))\overline{\dfrac{\partial c}{\partial z}(z)}}{\dfrac{\partial c}{\partial 
z}(z)} = \dfrac{\partial }{\partial \overline{z}} 
\left( \dfrac{\partial F}{\partial t}(0,z)\right)$
\end{lemma}
\proof 
To get the first equality, we begin by taking a derivative in $t$ in the expression of $\mu$ from the 
previous lemma to obtain
$$\dfrac{\partial \mu}{\partial t}(t,z)=
\dfrac{q \dfrac{\partial \eta}{\partial  \overline{w}}(c(z)) \overline{\dfrac{\partial c}{\partial 
z}(z)} }{\dfrac{\partial c}{\partial z}} \dfrac{\partial}{\partial t} 
\left( \dfrac{t}{1+tq\dfrac{\partial \eta}{\partial w}(c(z)) }  \right)=$$

$$= \dfrac{q \dfrac{\partial \eta}{\partial  \overline{w}}(c(z)) \overline{\dfrac{\partial 
c}{\partial z}(z)} }{\dfrac{\partial c}{\partial z}} \cdot
\dfrac{1+tq\dfrac{\partial \eta}{\partial w}(c(z))-tq\dfrac{\partial \eta}{\partial w}(c(z))}
{(1+tq\dfrac{\partial \eta}{\partial w}(c(z)))^2}=$$
  
$$= \dfrac{q \dfrac{\partial \eta}{\partial  \overline{w}}(c(z)) \overline{\dfrac{\partial 
c}{\partial z}(z)} }{\dfrac{\partial c}{\partial z}} \cdot
\dfrac{1}{\left(1+tq\dfrac{\partial \eta}{\partial w}(c(z))\right)^2}$$
and then we evaluate at $t=0$ to get the first identity. \par
For the second equality, let us recall the relation $c(z)+tq\eta(c(z))=H(t,c(z))=c_t(F(t,z))$. 
Differentiating in  $t$  we obtain
$$ q\eta(c(z))= \dfrac{\partial c_t}{\partial u}(F(t,z))  \dfrac{\partial F}{\partial t}(t,z)$$
from which we get
$$ \dfrac{\partial F}{\partial t}(t,z) = 
\dfrac{q\eta(c(z)) }{\dfrac{\partial c_t}{\partial u}(F(t,z)) }$$
By definition we have $F(0,z)=z$ and $c_0(F(0,z))=c(z)$, hence at $t=0$ we obtain
$$ \dfrac{\partial F}{\partial t}(0,z) = 
\dfrac{q\eta(c(z))}{\dfrac{\partial c}{\partial z}(z) }$$
Differentiating in $\overline z$ we obtain then
$$\dfrac{\partial }{\partial \overline{z}} \left( \dfrac{\partial F}{\partial t}(0,z)\right)=  
\dfrac{q}{\dfrac{\partial c}{\partial z}(z)}  \dfrac{\partial }{\partial \overline{z}} \left( 
\eta(c(z)) \right)= 
\dfrac{q\dfrac{\partial \eta}{\partial \overline{w}}(c(z))\overline{\dfrac{\partial c}{\partial 
z}(z)}}{\dfrac{\partial c}{\partial z}(z)}$$
which proves the second identity.
\endproof

\begin{remark}\label{beltramiasbarprimitive}
Notice that the first order deformation at $t=0$ has an elementary $\overline{\partial}$-primitive: namely by the above computations we get
$$ \dfrac{\partial \mu}{\partial t}(0,z)=\dfrac{\partial }{\partial \overline{z}}
 \left( \dfrac{q\eta(c(z))}{\dfrac{\partial c}{\partial z}(z)}  \right)=\dfrac{q}{\dfrac{\partial c}{\partial z}(z)}  
\dfrac{\partial }{\partial \overline{z}} ( \eta(c(z)) ) $$
where the second equality comes from the fact that $c$ is holomorphic.
This will turn out to be useful in the following computations.
\end{remark}

\subsection{Contraction with quadratic differentials}
The goal of this section is to provide a proof of the main theorem (Theorem \ref{thm_intro1} in the Introduction, see Theorem \ref{thm_holomoves_can} below). This involves a criterion to 
decide whether a movement of branch points induces a deformation of the underlying complex 
structures which is non-trivial at the first order. Beltrami differentials are specifically designed to 
measure whether a deformation of a complex structure is trivial or not, so in order to do this, we 
are going to determine whether the family of Beltrami differentials $\mu_t$ from the previous 
section is trivial at first order; equivalently, we regard the expression 
$\dot{\mu}_0(z)=\frac{\partial \mu}{\partial t}(0,z)$ from Lemma \ref{lem_beltrami_firstorder} as a 
Beltrami differential on its own, representing the first order behavior of the family $\mu_t$ at 
$t=0$.\par

It is a standard fact in Teichm\"uller theory (see \cite{IT92} for a reference) that given any 
complex structure $X \in \mathcal T(S)$ there is a natural pairing between Beltrami differentials 
and holomorphic quadratic differentials on $X$ given by

$$\langle \alpha,\mu \rangle = \int_X \alpha \mu = \int_X \alpha(z) \mu(z) dzd\overline z $$

where $\alpha=\alpha(z)dz^2$ is a local expression for a holomorphic quadratic differential 
$\alpha$ and $\mu=\mu(z) d\overline z \left(\frac{\partial}{\partial z} \right)^{-1}$ is a local 
expression for a Beltrami differential. Moreover a Beltrami differential $\mu$  is trivial (hence 
the associated deformation of complex structure is trivial) precisely when $\langle \alpha,\mu 
\rangle =0$ for all holomorphic quadratic differentials $\alpha \in H^0(X,K_X^2)$.\par

We are now going to compute an expression for the pairing of the first order approximation 
$\dot{\mu}_0$ with an arbitrary holomorphic quadratic differential. Let us recall the notation 
and setting from the previous section: $p$ is a branch point of order $m-1$, local projective charts 
can be taken to be of the form $c(z)=z^m+o(z^m)$, and let us introduce the notation $\frac{\partial 
c}{\partial z}(z)=z^{m-1}g(z)$, where $g$ is a holomorphic function  such that $g(0)=m$.
\begin{proposition}\label{contraction}
 Let $\alpha \in H^0(X,K^2_X)$ be a holomorphic quadratic differential on $X$, and let 
$\alpha=\alpha(z)dz^2$ be its expression in the coordinate $z$ on $\varphi(A)$. Then 
$$\langle \alpha,\dot{\mu}_0 \rangle=\dfrac{2\pi i q }{(m-2)!} \left. \dfrac{\partial^{m-2} }{\partial 
z^{m-2}}\dfrac{\alpha(z)}{g(z)}\right|_{z=0}$$

\end{proposition}
\proof
We begin with
$$\langle \alpha,\dot{\mu}_0 \rangle=\int_S \alpha\dot{\mu}_0 = \int_A \alpha\dot{\mu}_0 =$$
where we restrict the integral over $A$ since $\mu_t$ is compactly supported inside it for any 
$t\in [0,1]$. Then we can go in local coordinates in $\varphi(A)$
$$=\int_{\varphi(A)} \alpha(z) \dot{\mu}_0(z)dzd\bar{z}=
\int_{\varphi(A)} \alpha(z) \dfrac{q}{\dfrac{\partial c}{\partial z}(z)}  
\dfrac{\partial }{\partial \overline{z}} ( \eta(c(z)) ) dzd\bar{z}$$
where the last equality comes from the above Remark \ref{beltramiasbarprimitive}. Also notice that the second integral can actually be restricted to $\varphi(B\setminus C)$ because $\eta$ is constant on $c(\varphi(C))=d(C)$ and on $c(\varphi(A\setminus B))=d(A\setminus B)$.\par
We now observe that, since $\alpha$ is holomorphic, we have
$$d\left(\dfrac{q\alpha(z)}{\dfrac{\partial c}{\partial z}(z)}\eta(c(z))dz\right) 
=-\dfrac{q\alpha(z)}{\dfrac{\partial c}{\partial z}(z)} \dfrac{\partial }{\partial \overline{z}} ( 
\eta(c(z)) )dzd\bar{z} $$

Then we can resume the computation from above and obtain
$$\langle \alpha,\dot{\mu}_0 \rangle=\dots= -\int_{\varphi(B\setminus C)} d\left(\dfrac{q\alpha(z)}{\dfrac{\partial c}{\partial 
z}(z)}\eta(c(z))dz\right)  =  $$
to which we now apply Stokes Theorem to get
$$= -\int_{\varphi(\partial B)} \dfrac{q\alpha(z)}{\dfrac{\partial c}{\partial 
z}(z)}\eta(c(z))dz + \int_{\varphi(\partial C)} \dfrac{q\alpha(z)}{\dfrac{\partial c}{\partial 
z}(z)}\eta(c(z))dz = 
\int_{\varphi(\partial C)} \dfrac{q\alpha(z)}{z^{m-1}g(z)}dz$$
where the last equality comes from the fact that, by definition, $\eta=0$ in the first 
integral and $\eta=1$ in the second, and from the definition of $g(z)=\frac{\partial c}{\partial 
z}(z)z^{1-m}$. Now observe that everything inside the last integral is holomorphic, therefore we 
can apply Cauchy's integral formula to obtain the desired expression
$$\langle \alpha,\dot{\mu}_0 \rangle= \int_{\varphi(\partial C)} \dfrac{q\alpha(z)}{z^{m-1}g(z)}dz = \dfrac{2\pi i q }{(m-2)!} \left. \dfrac{\partial^{m-2} }{\partial 
z^{m-2}}\dfrac{\alpha(z)}{g(z)}\right|_{z=0}$$
\endproof

\begin{example}\label{ex_contraction_g=2}
 In the case of a simple branch point, i.e. $m=2$, we have $c(z)=z^2+o(z^2)$, hence $g(0)=2$ and 
the above formula reduces to 
 $$\langle \alpha,\dot{\mu}_0 \rangle=\pi i q \alpha(0)$$
\end{example}

We are now ready to prove the main result.
Recall from Definition \ref{def_holo_fibers} that if $k\in \mathbb N$ 
and $\lambda=(\lambda_1,\dots,\lambda_n)$ is a partition of $k$, then a structure $\sigma \in 
\mathcal{M}_{k,\rho}$ belongs to the $\lambda$-stratum $\mathcal{M}_{\lambda,\rho}$ if its branching 
divisor is of the form $div(\sigma)=\sum_{j=1}^n \lambda_jp_j$; in the above notations, a branch 
point of order $\lambda_j$ has a local chart of the form $z \mapsto 
z^{\lambda_j+1}+o(z^{\lambda_j+1})$.  \par

Since we are considering only deformations which preserve the structure of the branching divisor, 
i.e. do not leave the stratum $\mathcal{M}_{\lambda,\rho}$, we find it convenient to denote by 
$\pi^{\lambda}:\mathcal{M}_{\lambda,\rho}\to \mathcal{T}(S)$  the restriction of the Teichm\"uller map to the the 
$\lambda$-stratum $\mathcal{M}_{\lambda,\rho}$ of the holonomy fiber $\mathcal{M}_{k,\rho}$.

\begin{theorem}\label{thm_holomoves_can}
Let $\rho:\pi_1(S)\to\pslc$ be non-elementary, $k\leq 2g-2$, $\lambda=(\lambda_1,\dots,\lambda_n)$ be a partition of 
$k$; let $\sigma \in \mathcal{M}_{\lambda,\rho}$ and let $div(\sigma)=\sum_{j=1}^n \lambda_jp_j$ be its branching 
divisor.
If $\sigma$ is a critical point for $\pi^{\lambda}$, then $div(\sigma)$ is a canonical divisor 
on the underlying Riemann surface $X=\pi(\sigma)$.

\end{theorem}
\proof
As recalled at the beginning of this section, it is a classical fact that a 
deformation of a complex structure on a surface is trivial if and only if the contraction of its  
Beltrami differential with all holomorphic quadratic differentials is trivial. \par

Local deformations of $\sigma$ inside $\mathcal{M}_{\lambda,\rho}$ are precisely given by movements 
of branch points, according to Remark \ref{rmk_move_in_stratum}.
If $(A_j,d_j)$ is a local projective chart at $p_j$, for $j=1,\dots,n$, with $A_i\cap A_j=\varnothing$ for $i\neq j$, then we know from Remark
\ref{rmk_move_param} that these deformations are parametrized by the choice of a point  
$q=(q_1,\dots,q_n) \in \prod_{j=1}^n d_j(A_j)$. \par

We are now going to assume by contradiction that the divisor $div(\sigma)$ is not canonical on $X$, 
and show that all these deformations change the underlying complex structure, by showing that they 
are non-trivial at first order. In order to show this we will exploit the aforementioned natural 
pairing with holomorphic quadratic differentials, and the formula from Proposition \ref{contraction}.\par

If $(A_j,\varphi_j)$ are the local complex charts around $p_j$ associated to $(A_j,d_j)$ in the 
sense of Lemma \ref{lem_projchartincomplexchart}, with coordinate $z_j$, and if $c_j=d_j\circ 
\varphi_j^{-1}$ and 
$g_j(z_j)=\frac{\partial c_j}{\partial z_j}(z_j)z_j^{-\lambda_j}$, then by Proposition 
\ref{contraction} we can write the 
contraction between any holomorphic quadratic differential $\alpha \in H^0(X,K_X^2)$ and the first order approximation 
$\dot{\mu}_0$ of the Beltrami differential of this deformation as 

$$<\alpha,\dot{\mu}_0 > = \sum_{j=1}^n \dfrac{2\pi i q_j }{(\lambda_j-1)!} \left. 
\dfrac{\partial^{\lambda_j-1}  }{\partial z_j^{\lambda_j-1}} \dfrac{\alpha_j}{g_j} (z_j)\right|_{z_j=0}$$
where $\alpha_j$ is the local expression of $\alpha$ in the coordinate $(A_j,\varphi_j)$ at $p_j$.
Notice that for any $s\geq 0$ we have
$$\dfrac{\partial^{s}}{\partial z_j^{s}}\dfrac{\alpha_j}{g_j}(z_j)=
\dfrac{1}{g_j(z_j)} \left( \dfrac{\partial^{s} \alpha_j}{\partial z_j^{s}}(z_j) - 
\sum_{l=0}^{s-1}  {{s}\choose{l}} \dfrac{\partial^{l}}{\partial z_j^{l}}\dfrac{\alpha_j}{g_j}(z_j)
\dfrac{\partial^{s-l}g_j}{\partial z_j^{s-l}}(z_j) \right)$$
so that the vanishing of derivatives of $\alpha_j$ recursively implies the vanishing of derivatives of 
$\frac{\alpha_j}{g_j}$. Recalling that by definition $g_j(0)=\lambda_j+1$, in particular we obtain the following
$$\dfrac{\partial^{s} }{\partial z_j^{s}}\dfrac{\alpha_j}{g_j}(0)=
\left\lbrace \begin{array}{ll}
 \dfrac{1}{g_j(0)} \dfrac{\partial^{s} \alpha_j}{\partial z_j^{s}}(0)=\dfrac{a_js!}{\lambda_j+1}  & \textrm{if } 
ord_0(\alpha_j)=s \\
 0 & \textrm{if } ord_0(\alpha_j)>s\\ \end{array}\right.$$
for some constants $a_j\neq 0$.
If the movement of branch points is not trivial, then at least one of the points $q_j$ is not zero; let us assume for 
simplicity $q_1\neq 0$. We now claim that, under the assumption that $div(\sigma)$ is not canonical 
on $X$, it is possible to find a holomorphic quadratic 
differential $\alpha \in H^0(X,K_X^2)$ which has a zero of order exactly $\lambda_1-1$ at $p_1$ (not higher!) and a zero 
of order at least $\lambda_j$ at $p_j$ for $j=2,\dots,n$. For such a differential it is easily obtained from the above 
computations that
$$\dfrac{\partial^{\lambda_j-1} }{\partial z_j^{\lambda_j-1}}\dfrac{\alpha_j}{g_j}(0)=
\left\lbrace \begin{array}{ll}
 \dfrac{a_1(\lambda_1-1)!}{\lambda_1+1}  & j=1 \\
 0 & j=2,\dots,n\\ \end{array}\right.$$
We conclude that for such a differential we have $<\alpha,\dot{\mu}_0 >= \frac{2\pi i q_1 a_1}{\lambda_1+1}\neq 0$, 
which implies that deformation is non-trivial at first order, hence $\sigma$ can not be a critical 
point for $\pi^\lambda$, which is the desired contradiction. \par
To complete the proof let us prove the above claim about the existence of the required quadratic 
differential. If $E$ is any divisor on $X$, let us denote by $Q(X,E)=\{\alpha \in 
H^0(X,K_X^2) \ | \ div(\alpha)+E\geq 0\}$; when $E<0$ this is the space of holomorphic quadratic differentials vanishing 
along $E$ with at least the multiplicity prescribed by $E$. Let $D=-div(\sigma)$ and $D'=p_1+D$; notice $deg(D)=-k<0$ by definition. We 
have an obvious inclusion
$$Q(X,D) \subseteq Q(X,D')$$
An element of $Q(X,D)$ vanishes at $p_j$ with order at least $\lambda_j$ for all $j=1,\dots,n$, while an element of 
$Q(X,D')$ vanishes at $p_j$ with order at least $\lambda_j$ for all $j=2,\dots,n$, but only with order at least 
$\lambda_1-1$ at $p_1$. So an element of $Q(X,D_1)\setminus Q(X,D)$ is precisely a holomorphic 
quadratic differential 
which vanishes on $p_1$ with order exactly $\lambda_1-1$, and vanishes on $p_j$ with order at least $\lambda_j$ for 
$j=2,\dots,n$. The claim then follows if we prove that the above inclusion has positive codimension. The Riemann-Roch 
formula allows to compute that
\begin{itemize}
 \item if $div(\sigma)$ is not canonical, then $dim(Q(X,D))=3g-3-k$ and $dim(Q(X,D'))=3g-3-k+1$,
 \item if $div(\sigma)$ is canonical, then $dim(Q(X,D))=dim(Q(X,D'))$.
\end{itemize}
which concludes the proof of the claim.
\endproof

The remaining results of this section are straightforward consequences of the previous Theorem.

\begin{corollary}\label{cor_fibersarecan}
Let $k\leq 2g-2$ and let $X\in \mathcal T(S)$. If $F=\pi^{-1}(X)\subseteq \mathcal{M}_{k,\rho}$ has 
positive dimension, then for any $\sigma \in F$ we have that $div(\sigma)$ is canonical on $X$.
\end{corollary}
\proof
Let $\sigma \in F$.
First assume that a neighborhood $\Omega$ of $\sigma$ in $F$ is fully contained in a stratum 
$\mathcal M_{\lambda, \rho}$, for some partition $\lambda$ of $k$. Then $\pi^\lambda( \Omega)=X$, in 
particular $\sigma$ is critical for $\pi^\lambda$. By Theorem \ref{thm_holomoves_can} we are 
done.\par
So now assume that any neighborhood of $\sigma$ in $F$ contains structures from 
higher-dimensional strata; then consider an arbitrarily small deformation of $\sigma$ that splits 
branch points on $\sigma$ to jump into the top-dimensional stratum of the induced stratification on 
$F$ (see Remark \ref{rmk_move_in_stratum}). Considering the structures obtained in this way, we can 
argue as before, and apply Theorem \ref{thm_holomoves_can} to show that they are all critical for 
the projection from their stratum, hence they are all  canonically branched. This shows 
that $\sigma$ has full connected neighborhood $\Omega$ in $F$  in which canonically branched structures 
fill an open dense subspace. Since being canonical is a closed condition on divisors, this implies 
that $\sigma$ is canonically branched too (as well as every other structure in $\Omega$).
\endproof
\begin{example}
 A major motivating example is given by the 1-dimensional submanifold $\Sigma_{(X,A)}\subseteq F= 
\pi^{-1}(X)$ where $(X,A)$ is an $\mathfrak{sl}_2$-system (see Section \ref{sec_fromODEtoBPS} for 
details about the construction). Structures on $ \Sigma_{(X,A)}$ have the same underlying complex 
structure and are generically contained in the same stratum (see Remark 
\ref{rmk_fiber_transverse_stratum} for more subtle aspects), and they have already been shown to be 
branched on a canonical divisor in Proposition \ref{prop_ODEBPS_can}.
\end{example}

\begin{remark}\label{rmk_fiber_transverse_stratum}
It is possible for positive-dimensional fibers and positive-dimensional strata to 
intersect at isolated points (which is covered by the second case of the previous proof). 
For instance if $(X,A)$ is a $\mathfrak{sl}_2$-system
 of genus $g=2$ (see Section \ref{sec_fromODEtoBPS}), then the induced rational curve 
$\Sigma_{(X,A)}$ is a 1-dimensional fiber over $X$ of the 2-dimensional complex manifold $\mathcal 
M_{2,\rho}$; it is  generically  contained in the principal stratum $\mathcal M_{(1,1),\rho}$ 
but it may intersect the minimal stratum $\mathcal M_{(2),\rho}$ in six isolated points (corresponding 
to structures with a double branch point at one of the six Weierstrass points of $X$).
\end{remark}

The following is based on an observation present in \cite[Corollary 6.3]{CDHL16}, where it is shown 
that if a fiber of $\pi:\mathcal{M}_{2g-2,\rho}\to \mathcal T(S)$ contains at least three points, 
then it is actually   the rational curve of some $\mathfrak{sl}_2$-system (see Section 
\ref{sec_fromODEtoBPS}).
\begin{corollary}\label{cor_compactsubarecan}
Let $k\leq 2g-2$. If $Z\subseteq \mathcal{M}_{k,\rho}$ is a compact complex submanifold of positive 
dimension, then  for any $\sigma \in Z$ we have that $div(\sigma)$ is canonical on the underlying 
Riemann surface.
\end{corollary}
\proof
Since Teichm\"uller space is a Stein manifold, $Z$ is compact and $\pi:\mathcal{M}_{k,\rho}\to 
\mathcal{T}(S)$ is holomorphic, we have that the restriction of $\pi$ to $Z$ must be constant. 
In other words $Z$ must be contained in a fiber $\pi^{-1}(X)$ for some $X \in \mathcal{T}(S)$.
Then we can apply Corollary \ref{cor_fibersarecan} to conclude.
\endproof

\subsection{The sub-canonical range}\label{sec_subcan}
In this section we apply results from the previous section to show that the  Teichm\"uller  map
$\pi:\mathcal{M}_{k,\rho} \to \mathcal T(S)$ is quite well-behaved when the branching order 
is sub-canonical, i.e. $k<2g-2$.\par

\begin{corollary}\label{cor_subcan}
Let $\rho:\pi_1(S)\to\pslc$ be non-elementary, $k< 2g-2$ and $\lambda$ a partition of $k$. Then the following hold:
\begin{enumerate}
 \item $\pi^\lambda:\mathcal{M}_{\lambda,\rho}\to \mathcal{T}(S)$ is an injective immersion;
 \item $\mathcal{M}_{k,\rho}$ contains no positive-dimensional compact complex submanifolds.
\end{enumerate}
\end{corollary}
\proof
${}$
\begin{enumerate}
 \item The injectivity of $\pi^\lambda$ follows from the global injectivity of $\pi$, which is a classical result of 
Hejhal (see \cite[Theorem 15]{HE75}, which appears as Theorem \ref{thm_hejhal} above). 
Then assume by contradiction that $\pi^\lambda$ has a critical point $\sigma$. By Theorem 
\ref{thm_holomoves_can} we obtain that $\sigma$ must be canonically branched. But this is absurd, 
simply because the branching order of $\sigma$ is less than the degree of a canonical divisor. 
\item Analogously, by Corollary \ref{cor_compactsubarecan} structures on a compact complex 
submanifold should be canonically branched, which once again is impossible just because $k<2g-2$.
\end{enumerate}
\endproof

In light of the previous statement, it would then  be interesting to consider the following questions:
\begin{question}
Is $\pi$ immersive?
\end{question}
\begin{question}
Are $\pi$ and/or $\pi^\lambda$ proper?
\end{question}

Affirmative answers would imply (in view of Corollary \ref{cor_subcan}) that when $k<2g-2$ the 
manifolds $\mathcal{M}_{k,\rho}$ (or at least their strata $\mathcal{M}_{\lambda,\rho}$) embed as a closed smooth complex submanifolds of Teichm\"uller space, hence are themselves Stein manifolds; point (2) in Corollary \ref{cor_subcan} suggests that no obvious obstruction to this arises from the geometry of submanifolds.

\section{Hyperelliptic structures}\label{sec_hyperelliptic}
The purpose of this last section is  to obtain a partial converse to the main result (Theorem \ref{thm_holomoves_can}) in the case of BPSs endowed with a special type of automorphism. We will provide an explicit construction of structures of this type and study natural deformations thereof.

\subsection{Projective hyperelliptic involutions}
Recall from Definition \ref{def_BPS_stuff}  that, given a BPS $\sigma 
\in \mathcal {BP} (S)$, a diffeomorphism $f$ of $S$ is said to be projective with respect to 
$\sigma$ if it is given by M\"obius transformations in local projective charts defining $\sigma$. 
Equivalently, if $f$ acts on developing maps for $\sigma$ by post-composition with some M\"obius 
transformation in the following sense (see \cite[\S 2]{FR19} for details): let $dev$ be a developing 
map for $\sigma$, then $f$ is projective for $\sigma$ if and only if $\exists \ g \in \pslc$ and a 
lift $\widetilde f$ of $f$ to the universal cover such that $dev\circ \widetilde{f}^{-1}=g\circ 
dev$; in particular $dev\circ \widetilde{f}^{-1}$ is another developing map for $\sigma$.

Also recall from Remark \ref{rmk_teich} that a BPS has an underlying complex structure, and notice 
that projective diffeomorphisms are in particular biholomorphic with respect to it. The following 
lemmas give some elementary properties of projective diffeomorphisms.

\begin{lemma}\label{lem_permute_brpts}
Let $\sigma\in \mathcal {BP} (S)$ and let $f:S\to S$ be a projective diffeomorphism. Then $f$ 
permutes the branch points of $\sigma$, preserving their order. 
\end{lemma}
\proof
Let $dev_1$ be a developing map for $\sigma$, and let $dev_2=dev_1\circ \widetilde{f}^{-1}=g\circ dev_1$ for some 
lift of $f$ to the universal cover and some $g\in \pslc$. For any point $p \in \widetilde{S}$ in the universal cover  we 
have
$$\partial_z dev_1(\widetilde{f}^{-1}(p)) \partial_z \widetilde{f}^{-1}(p)=\partial_z 
dev_2(p)=\partial_w g (dev_1(p)) \partial_z dev_1(p)$$
for complex coordinates $z$ on $\widetilde S$ and $w$ on $\cp$.
Since $\widetilde{f}$ and $g$ are biholomorphisms, this means that $p$ is a branch point for $dev_1$ 
if and only if it is a branch point for $dev_2$, if and only if $\widetilde{f}^{-1}(p)$ is a branch 
point for $dev_1$. Branching orders are precisely given by the order of ramification of developing 
maps, so the same equation shows that they must be preserved.
\endproof

We focus here on structures which admit a special type of projective symmetry. Recall that a Riemann 
surface $X$ of genus $g$ is hyperelliptic if it admits a holomorphic involution $J$ with $2g+2$ fix 
points, known as the Weierstrass point of $X$; notice that if follows from the Riemann-Hurwitz formula that $X/J\cong \cp$.
\begin{definition}\label{def_hyperellipticBPS}
A BPS is said to be hyperelliptic if it admits a projective involution $J$ with $2g+2$ fix points. 
The involution will be called a projective hyperelliptic involution.
A hyperelliptic BPS is non-degenerate if the projective hyperelliptic involution $J$ does not fix 
any of the branch points. 
If $\mathcal X \subseteq \mathcal {BP}(S)$ is any subspace, then we will denote by 
$\mathcal {HX}$ (respectively $\mathcal {HX}^\sharp$) the locus of hyperelliptic 
(respectively  non-degenerate hyperelliptic) structures in $\mathcal {X}$.
\end{definition}
\begin{example}
 Any hyperelliptic hyperbolic surface can be seen as a hyperelliptic BPS without branch points: 
the hyperelliptic involution $J$ is isometric, hence projective.
\end{example}
\begin{example}
For a branched example, consider the trivial BPS induced on a hyperelliptic Riemann surface $X$ by 
the quotient map $X\to X/J=\cp$; this structure has $2g+2$ simple branch points and trivial holonomy. For a less trivial example, let $\sigma_0$ be a genus zero Euclidean 
orbifold with four points of angle $\pi$, and take a double cover $\sigma \to \sigma_0$ branched at 
the four cones and also at two smooth points; the result is a hyperelliptic BPS on a surface of genus 
2. An analogous construction works in every genus, but notice that all these structures are 
degenerate, as that they are branched precisely at Weierstrass points. For examples of 
non-degenerate hyperelliptic BPSs we refer to Section \ref{sec_hyperelliptic_bubbling} below.
\end{example}

We collect here for future reference some facts about hyperelliptic structures that just follow from the definition and 
Lemma \ref{lem_permute_brpts}.
\begin{lemma}\label{lem_hyperelliptic_stuff}
Let $\sigma \in \mathcal{HBP}(S)$, with projective hyperelliptic involution $J$, and underlying 
complex structure $X$. Then
\begin{enumerate}
 \item $X$ is a hyperelliptic Riemann surface with hyperelliptic involution $J$.
 \item $\sigma \in \mathcal{HBP}(S)^\sharp$   if and only if none of the 
branch points of $\sigma$ is a Weierstrass  point for $X$.
\item If $\sigma \in \mathcal{HBP}(S)$, $p $ is a branch point, and $(A,d)$ is a projective 
chart at it with $d(p)=0$, then $p'=J(p)$ is a branch point of the same order, and
$(A',d')=(J(A),d\circ J)$ is a projective chart at it with $d'(p')=0$ and $d'(A')=d(A)$.
 
\end{enumerate}
\end{lemma}

For non-degenerate hyperelliptic BPS we can introduce preferred neighborhoods inside the 
holonomy fibers as follows.
Let $\sigma \in \mathcal{HM}^\sharp_{\lambda,\rho}\subset \mathcal M_{k,\rho}$ for some partition $\lambda$ of $k$. Thanks to Lemma \ref{lem_hyperelliptic_stuff} 
the branch points of $\sigma$ are in even number and of the form $p_1,J(p_1),\dots,p_n,J(p_n)$, where $p_i\neq 
J(p_i)$, and $n$ is half of the length of the partition $\lambda$; in particular branch points are 
simple precisely when $2n=k$, but we are not requiring this.

\begin{definition}\label{def_stdnbd}
A standard neighborhood of $\sigma \in \mathcal{HM}^\sharp_{\lambda,\rho}$ is the neighborhood of 
$\sigma$  inside $\mathcal{M}_{\lambda,\rho}$ which is defined, in the notation of Definition 
\ref{def_move},  by
$$\Omega(\sigma)=\{Move(\sigma, z_1,w_1,\dots, z_n,w_n) \ | \ z_i,w_i \in d_i(A_i), i=1,\dots,n \} 
$$ 
where $(A_i,d_i)$ is a projective chart at $p_i$, for $i=1,\dots,n$.
We also define the diagonal slice of it to be the half-dimensional subspace given by
$$\Delta\Omega(\sigma) = \{ Move(\sigma,z_1,w_1,\dots,z_n,w_n ) \in \Omega(\sigma) \  | \ z_i=w_i, i=1,\dots,n 
\} $$ 
\end{definition}
Naively speaking, this is the neighborhood of structures obtained by moving branch points in such 
a way that paired branch points $p_i,J(p_i)$ remain inside the common chart $d_i(A_i)$, and the 
diagonal slice is the subspace obtained by moving them ``consistently with respect to $J$'' (i.e. 
$J$-equivariantly).
Notice that 
$$\Omega(\sigma)\cong \prod_{i=1}^n (d_i(A_i)\times d_i(A_i))  \ \ , \ \ \Delta\Omega(\sigma) \cong 
\prod_{i=1}^n d_i(A_i)$$
Of course   a  standard neighborhood depends on a specific choice of projective charts at each 
branch point $p_i$ for $i=1,\dots,n$; but in any case all the possible resulting standard 
neighborhoods provide a fundamental system of neighborhoods of $\sigma$  for the topology of 
$\mathcal{M}_{\lambda,\rho}$ (compare Remarks \ref{rmk_move_param} and \ref{rmk_move_in_stratum}). 
Hence we suppress the dependence in our notation, and whenever talking about a standard 
neighborhood or diagonal slice we will implicitly assume that such a choice of local charts has been made.

We conclude this section by showing that diagonal slices are actually contained in the non-degenerate
hyperelliptic locus. This is illustrated by Example \ref{ex_hypbub} and Remark \ref{rmk_hypbub_deformation} below. 

\begin{proposition}\label{prop_hyp_slice}
If $\sigma \in \mathcal{HM}^\sharp _{\lambda,\rho}$,
then  $\Delta\Omega(\sigma)\subseteq \mathcal{HM}^\sharp_{\lambda,\rho}$.
\end{proposition}
\proof
Recall from Section \ref{sec_surgeries} that local deformations of $\sigma$ in $\mathcal{M} _{\lambda,\rho}$ 
are given by moving branch points, and that this can be obtained by postcomposing charts  around branch points with a  compactly supported isotopies. By Lemma \ref{lem_hyperelliptic_stuff} we can choose charts to be of the form  $(A_i,d_i)$ at $p_i$ and  $(J(A_i),d_i\circ J)$ at $J(p_i)$; let $\Omega(\sigma)$ be the induced standard neighborhood. Then deformations along the diagonal slice $\Delta \Omega(\sigma)$ correspond to using the same isotopy $H_i$ of $d_i(A_i)$ for both $p_i$ and $J(p_i)$.
The hyperelliptic involution of $\sigma$, seen just as a diffeomorphism of $S$, is given by the identity in these preferred charts throughout the deformation; in particular it is projective for all of them. This proves that all the structures along the diagonal slice are hyperelliptic. Moreover they are non-degenerate because the involution never fixes any of the branch points by construction.
\endproof

\begin{remark}
It has been observed in \cite[Lemma 10.1.(iii)]{CDHL16} that in the genus 2 case for any non-elementary representation  
$\rho$ the hyperelliptic involution induces a holomorphic involution of $\mathcal M_{2,\rho}$ without isolated fixed 
points. Proposition \ref{prop_hyp_slice} provides a generalization to higher genus and higher number of branch points: 
when the holonomy fiber $\mathcal M_{k,\rho}$  supports an action by a hyperelliptic involution, the hyperelliptic 
locus is  at least half-dimensional in every stratum.
 
\end{remark}

\subsection{Hyperelliptic bubblings} \label{sec_hyperelliptic_bubbling}
Here we present a general construction to produce examples of non-degenerate hyperelliptic BPSs 
(see Definition \ref{def_hyperellipticBPS}). 

\begin{example}\label{ex_hypbub}
Let $\sigma \in \mathcal{HM}_{0,\rho}$, with involution $J$, and let $\beta:[0,1]\to \sigma$ be a 
bubbleable arc on $\sigma$ (see Definition \ref{def_bubbling}). Assume that $J$ maps $\beta$ to 
itself swapping its two endpoints. Then bubbling $\sigma$ along $\beta$ results in a non-degenerate 
hyperelliptic structure $Bub(\sigma, \beta) \in \mathcal{HM}^\sharp_{2,\rho}$.
Indeed  let $dev$ be a developing map and let $g\in \pslc$ be an order 2 M\"obius transformation realizing $J$ through $dev$, in the sense that $dev\circ \widetilde{J}^{-1}=g\circ dev$ where $\widetilde J$ is a lift of $J$ to the universal cover.
Since $\beta$ is $J$-invariant, its developed image $\widehat{\beta}$ is $g$-invariant, therefore the bubble 
$\cp\setminus \widehat{\beta}$ carries a natural projective involution given by the restriction of 
$g$. We first restrict $J$ to an involution of $\sigma \setminus \beta$, and then just extend this to the whole $Bub(\sigma, \beta)$ by using $g$ inside the bubble.
\end{example}

\begin{remark}\label{rmk_hypbub_deformation}
This construction can be deformed to give rise to a 1-dimensional family inside 
$\mathcal{HM}^\sharp_{2,\rho}$: for instance just extend $\beta$ as desired at one endpoint and 
extend $\beta$ at the other endpoint in the only way that produces a $J$-invariant arc; then 
perform a bubbling on the new arc. 
In any standard neighborhood for $\sigma$ this family is given by the diagonal slice $\Delta 
\Omega(\sigma)$ (see Definition \ref{def_stdnbd} and Proposition \ref{prop_hyp_slice}).
\end{remark}

\begin{remark}
This construction applies every time some  $\sigma \in \mathcal{HM}_{0,\rho}$ is available, for 
instance
\begin{itemize}
 \item when $\rho$ is a hyperelliptic Fuchsian representation (i.e. such that $\mathbb H^2 / \rho$ 
is a hyperelliptic hyperbolic surface). In this case just take $\sigma = \mathbb H^2 / \rho$. Representations of this type exist in every genus $g\geq 2$.
\item in genus $g=2$ every unbranched projective structure is hyperelliptic (see for instance 
\cite[Example 4.11]{FR19}); so for any non-elementary liftable representation $\rho$ we can choose 
any $\sigma\in \mathcal{M}_{0,\rho}=\mathcal{HM}_{0,\rho}$.
\end{itemize}
Of course this construction can also be iterated to produce examples with more branch points, 
possibly with higher branching order too.
\end{remark}

\subsection{Critical movements for hyperelliptic structures}
In this section we prove a partial converse to Theorem  \ref{thm_holomoves_can} under the assumption that the structure 
is hyperelliptic. In particular we will show that the critical direction is always 1-dimensional, and provide an 
explicit equation for this line in suitable standard neighborhood coordinates (introduced in Definition
\ref{def_stdnbd}). Recall from that section that if $\lambda=(\lambda_1,\dots,\lambda_n)$ is a partition of $k$, 
$\sigma \in \mathcal{HM}^\sharp_{\lambda,\rho}$  and $\Omega(\sigma)$ is a standard neighborhood of $\sigma$, then the 
diagonal slice of $\Omega(\sigma)$ is defined to be

$$\Delta\Omega(\sigma) = \{ Move(\sigma,z_1,w_1,\dots,z_n,w_n ) \in \Omega(\sigma) \  | \ z_i-w_i=0, i=1,\dots,n 
\} $$

we also introduce a natural complement, which we call the opposite diagonal slice, defined as

$$\Delta^\perp \Omega(\sigma) = \{ Move(\sigma,z_1,w_1,\dots,z_n,w_n ) \in \Omega(\sigma) \  | \ z_i+w_i=0, i=1,\dots,n 
\} $$

In the rest of the section we will consider the case $k=2g-2$, and the generic case of structures contained in the 
principal stratum, corresponding to the partition $\lambda=(1,\dots,1)$, i.e. $n=g-1$ and all branch points are 
simple branch points (order 1). Moreover $\sigma \in  \mathcal{HM}_{(1,\dots,1),\rho}$ will be a hyperelliptic 
structure in the principal stratum, $J$ the hyperelliptic involution and $X=\pi(\sigma)$ the underlying Riemann 
surface. We know that the branching divisor will have the form

$$div(\sigma)=p_1+J(p_1)+\dots+p_{g-1}+J(p_{g-1})$$

Moreover the following is true.

\begin{lemma}\label{lem_hyperelliptic_canonical_divisor}
$div(\sigma)$ is a canonical divisor on $X$ if and only if $\sigma \in  \mathcal{HM}^\sharp_{(1,\dots,1),\rho}$.
\end{lemma}
\proof
On a hyperelliptic Riemann surface $X$ of genus $g$, the canonical divisor is linearly equivalent to the sum of $g-1$ 
fibers of the canonical map $X\to X/J=\cp$. Such a fiber can be given by a double Weierstrass point, or by the sum of 
two distinct non-Weierstrass points exchanged by the hyperelliptic involution. Since all the coefficients in our 
divisor are 1 (we are in the principal stratum), this divisor is canonical if and only if none of the points in it is a 
Weierstrass point.  The statement  then follows from Lemma \ref{lem_hyperelliptic_stuff} above.
\endproof

In the next proof we will need some standard facts about quadratic differentials on hyperelliptic Riemann surfaces. Let 
us fix an affine model for $X$, i.e. a realization of $X$ as the smooth compactification of  the plane affine curve 
$$C=\left\lbrace (x,y) \in \mathbb C^2 \ \left| \ y^2=P(x)=\prod_{j=1}^{2g+2}(x-w_k) \right. \right \rbrace$$
where $P\in \mathbb C[x]$ is a polynomial with distinct simple roots $w_j,j=1,\dots, 2g+2$. The points $W_j=(w_j,0)\in 
C$ are precisely the Weierstrass points,   the hyperelliptic involution is given by $J(x,y)=(x,-y)$, and $x$ is the 
quotient map. Notice that, since $P$ has even degree, $X$ has two points over $x=\infty$, and none of them is a Weierstrass 
point.\par

The hyperelliptic involution acts as a linear involution on the space of holomorphic quadratic differentials on $X$, 
with eigenvalues $\pm 1$, and splits it into the sum of the two eigenspaces $V_\pm$. Using the affine model they can be 
presented as follows:
$$H^0(X,K_X^2)=V_+\oplus V_-=\mathbb C[x]_{\leq 2g-2} \cdot \dfrac{dx^2}{y^2} \oplus \mathbb C[x]_{\leq g-3} 
\cdot \dfrac{dx^2}{y}$$
where $\mathbb C[x]_{\leq d} $ denotes the space of polynomials of degree 
at most $d$, and $\frac{dx^2}{y^2}$ (respectively $\frac{dx^2}{y}$) is a $J$-invariant (respectively 
$J$-anti-invariant) holomorphic differential on $X$. Notice $V_\pm$ have dimension respectively $2g-1$ and $g-2$; in 
particular all holomorphic quadratic differentials are $J$-invariant in genus $g=2$.

\begin{theorem}\label{thm_hyp_converse}
Let $\rho:\pi_1(S)\to \pslc$ be non-elementary, $\sigma \in  \mathcal{HM}_{(1,\dots,1),\rho} \subseteq \mathcal 
M_{2g-2,\rho}$ and let  $X=\pi(\sigma)$ be the underlying Riemann surface.
If  $div(\sigma)$ is a canonical divisor on $X$, then $\sigma$ is 
a critical point for $\pi=\pi^{(1,\dots,1)}$, with 1-dimensional critical direction.\par
More precisely, let $div(\sigma)=\sum_{j=1}^{g-1}(p_j+J(p_j))$ and $(A_i,d_i)$ be charts at $p_i$; then there exist non-zero 
constants $Q_1,\dots,Q_{g-2}\in \mathbb C^*$ such that the subspace of the induced standard neighborhood 
$\Omega(\sigma)$ defined by
 
$$ K\Omega(\sigma) = \left\lbrace Move(\sigma,z_1,w_1,\dots,z_{g-1},w_{g-1} ) \in \Omega(\sigma)    \left| 
\begin{array}{ll}
 z_i+w_i=0, & 1\leq i \leq g-1 \\
 z_j=Q_jz_{g-1}, & 1\leq j \leq g-2\\
\end{array}  \right. \right \rbrace$$
satisfies the following properties
\begin{enumerate}
 \item $K\Omega(\sigma)\subseteq \Delta^\perp\Omega(\sigma)$ (it is understood that $K\Omega(\sigma)= 
\Delta^\perp\Omega(\sigma)$ for $g=2$).
 \item $dim(K\Omega(\sigma))=1$.
 \item $T_\sigma K\Omega(\sigma) = \ker (d\pi )_\sigma$.
\end{enumerate}
\end{theorem}

In other words $K\Omega(\sigma)$ is precisely the critical locus, it is a line in standard neighborhood coordinates, 
and it is transversal to the diagonal slices.

\proof
To begin with we observe that $\pi=\pi^{(1,\dots,1)}$ at $\sigma$, just because $\sigma$ is assumed to belong to the 
principal stratum, which is an open dense submanifold of $\mathcal M_{2g-2,\rho}$. Moreover we know from Lemma 
\ref{lem_hyperelliptic_canonical_divisor} that $\sigma$ is non degenerate. So it makes sense to talk about 
standard neighborhoods. Now (1) is completely obvious from the definition of $\Delta^\perp \Omega(\sigma)$. To get (2) 
one can just observe that in these coordinates $K\Omega(\sigma)$ is given by $g-1+g-2=2g-3$ linear equations in $2g-2$ 
variables; since $Q_j\neq 0$ one sees that they are all independent. Indeed (1) holds for any choice of $Q_j$, and 
(2) holds for any choice of $Q_j\neq 0$. Therefore the theorem will be proven if we prove the existence of $Q_j\neq 0$ 
for which (3) is satisfied.\par

Let $(A_j,\varphi_j)$ be the local complex chart associated to $(A_j,d_j)$ in the sense of Lemma
\ref{lem_projchartincomplexchart}. Notice that $(J(A_j),d_j\circ J)$ and $(J(A_j),\varphi_j \circ J)$ are then local 
projective and complex chart at $J(p_j)$ respectively (see Lemma \ref{lem_hyperelliptic_stuff}), adapted to each other 
in the sense of Lemma \ref{lem_projchartincomplexchart}. \par

Let us pick a point $(z_1,w_1,\dots,z_{g-1},w_{g-1})\in \prod _{j=1}^{g-1} d_j(A_j)\times d_j(A_j)$, let $\mu_t$ be 
the Beltrami differential of the movement of branch point defined by
$$[0,1] \ni t\mapsto Move(\sigma, tz_1,tw_1,\dots,tz_{g-1},tw_{g-1})\in \Omega(\sigma)$$ 
and let $\dot{\mu}_0$ its derivative at $t=0$. As we have already done in Section \ref{sec_complexanalytic}, to test triviality of the Beltrami differential we are going to contract it with holomorphic quadratic differentials. We want to prove that movements along $K\Omega(\sigma)$ are precisely those for which the contraction with all holomorphic quadratic differentials vanishes; we are going to exploit the decomposition into invariant and anti-invariant differentials described above.\par
By Proposition \ref{contraction} (notice that all branch points are simple, and see also Example 
\ref{ex_contraction_g=2}) we can compute that for all $\alpha \in H^0(X,K_X^2)$ we have
$$<\dot{\mu}_0,\alpha>=\pi i \sum_{j=1}^{g-1} (\alpha_j(0)z_j+\alpha'_j(0)w_j) $$
where $\alpha_j(0)$ and $\alpha'_j(0)$ denote respectively the values of $\alpha$ at $p_j$ and $J(p_j)$ in the chosen 
complex coordinate (the ones adapted to the projective structure in the sense of Lemma \ref{lem_projchartincomplexchart}).
Notice that, since these coordinates are compatible with the action of $J$, if $\alpha \in V_\pm$ then $\alpha'_j(0)=\pm \alpha_j(0)$, for $j=1,\dots,g-1$. \par
In particular for any $\alpha \in V_+$ we have
$$<\dot{\mu}_0,\alpha> =  \pi i \sum_{j=1}^{g-1} \alpha_j(0)(z_j+w_j)  $$

We claim (*) that this vanishes $\forall \ \alpha \in V_+$ if and only if 
$$Move(\sigma, z_1,w_1,\dots,z_{g-1},w_{g-1})\in \Delta^\perp\Omega(\sigma)$$
To see this, notice that if $Move(\sigma, 
z_1,w_1,\dots,z_{g-1},w_{g-1})\in \Delta^\perp\Omega(\sigma)$ then $z_j+w_j=0$ for all $j=1,\dots,g-1$, so the sum clearly vanishes.
Vice versa if  $<\dot{\mu}_0,\alpha>=0$ for all  $\alpha \in V_+$, then in particular it has to vanish for the holomorphic quadratic differentials 
$\alpha_l$ given in the affine model by
$$\alpha_l= \prod_{k=1,k\neq l}^{g-1}(x-x_k)\dfrac{dx^2}{y^2} \ , \ l=1,\dots,g-1$$
where $p_j=(x_j,\sqrt{P(x_j)})$ and $J(p_j)=(x_j,-\sqrt{P(x_j)})$ are the coordinates of the branch points of $\sigma$ in 
the affine model for $X$. Notice that $\alpha_l$ vanishes at $p_j$ and $J(p_j)$ if and only if $j\neq l$; even if a priori the complex charts chosen above and the complex charts coming from the affine model are not the same, nevertheless the vanishing of a differential does not depend on the chosen chart, so we can safely conclude that $(\alpha_l)_j(0)=0$ if and only if $j\neq l$.
Then we can compute that
$$0=<\dot{\mu}_0,\alpha_l> = \pi i (\alpha_l)_l(0) (z_l+w_l)$$
which forces $z_l+w_l=0$ for $l=1,\dots,g-1$, and proves the claim (*).\par

If $g=2$, then every quadratic differential is invariant, hence we are done. If the genus is higher we need to take 
care of anti-invariant differentials.
Exploiting the affine model once again, let us consider the following differentials
$$\beta_l = \prod_{k=1,k\neq l}^{g-2}(x-x_k)  \dfrac{dx^2}{y} \ , \ l=1,\dots,g-2$$
If $g=3$, then it is understood that there is only one differential $\beta_1=\frac{dx^2}{y}$.
In any case $\{\beta_l \ | \ l=1,\dots,g-2\}$ form a basis for $V_-$ (see Remark \ref{rmk_basis_weird} at the end of this proof).
Moreover we have

$$( \beta_l)_j(0) = \left\lbrace 
\begin{array}{ll}
\dfrac{1}{\sqrt{P(x_j)}}\prod_{k=1,k\neq l}^{g-2}(x_{j}-x_k) =0 , & 1\leq j\leq g-2,j\neq l \\
\dfrac{1}{\sqrt{P(x_l)}}\prod_{k=1,k\neq l}^{g-2}(x_{l}-x_k) \neq 0 , & 1\leq j\leq g-2,j= l \\
\dfrac{1}{\sqrt{P(x_{g-1})}}\prod_{k=1,k\neq l}^{g-2} (x_{g-1}-x_k) \neq 0 , & j=g-1 \\
\end{array} \right. $$
where, once again, we only care about the vanishing pattern, as the actual non-zero values might be different in the chosen charts and in the chart coming from the affine model. Moreover by anti-invariance of $\beta_l$ we get that $$( \beta_l')_j(0) = - ( \beta_l)_j(0)$$

Let us consider a deformation $Move(\sigma, z_1,w_1,\dots,z_{g-1},w_{g-1})\in \Delta^\perp\Omega(\sigma)$, and let $\dot \mu_0$ be the induced Beltrami differential as above. Then we claim (**) that $<\dot{\mu}_0,\beta_l>=0 $ for $l=1,\dots,g-2$ if and only if  $Move(\sigma, z_1,w_1,\dots,z_{g-1},w_{g-1})\in K\Omega(\sigma)$, for a suitable choice of the constants defining $K\Omega(\sigma)$, which we will now compute. Indeed the contraction with these differentials can be computed to be

$$<\dot{\mu}_0,\beta_l> = 
\pi i ( ( \beta_l)_l(0)z_l  + ( \beta_l')_l(0) w_l  + ( \beta_{g-1})_{g-1}(0)z_{g-1}  + ( \beta_{g-1}')_l(0) w_{g-1}    ) =$$
$$= 2\pi i ( ( \beta_l)_l(0)z_l + ( \beta_{g-1})_{g-1}(0)z_{g-1}  ) $$
where the last equality follows from anti-invariance of $\beta_l$ and the fact that the deformation is along $\Delta^\perp\Omega(\sigma)$. Then one sees that all these quantities vanish precisely when the following system of linear equations is satisfied
$$  ( \beta_l)_l(0)z_l + ( \beta_{g-1})_{g-1}(0)z_{g-1} = 0 \ , \ l=1,\dots, g-2$$
Since all the coefficients here are non zero, setting $Q_l=-\dfrac{ ( \beta_{g-1})_{g-1}(0)}{( \beta_l)_l(0)}$ provides the constants required to prove the claim (**). The Theorem follows then, with this choice of constants, from claims (*) and (**), and the fact that invariant differentials together with $\{\beta_l \ | \ l=1,\dots,g-2\}$ generate the whole space of holomorphic quadratic differentials on $X$.
\endproof

\begin{remark}\label{rmk_basis_weird}
Notice that $\beta_l \in V_-$ because the polynomial defining it has degree $g-3$, since in the product the index is bounded to be $k\leq g-2$. This might look a bit odd, as it introduces an asymmetry in the role of the variables, but allowing $k\leq g-1$ would result in a polynomial of degree $g-2$, which would give rise to a meromorphic differential with poles at infinity in the affine model. It should be clear from the statement of the Theorem, in particular from the equations defining $K\Omega(\sigma)$ that this asymmetry is not substantial, and the variable $z_{g-1}$ which is made special by the choice of these differentials turns out to be just a parameter for $K\Omega(\sigma)$, and could be replaced by any of the other variables.
\end{remark}

\begin{remark}
As already observed in Remark \ref{rmk_fiber_transverse_stratum} it is possible for a hyperelliptic 
BPS in the stratum $\mathcal{HM}_{(2),\rho}$ to be canonically branched without being a critical 
point for the projection $\pi^{(2)}$; in other words the analogous statement for the minimal stratum fails.
\end{remark}

\begin{remark}
In particular a diagonal slice $\Delta\Omega(\sigma)$ at a non-degenerate hyperelliptic BPS 
$\sigma \in \mathcal{M}_{(1,\dots,1),\rho}\subseteq \mathcal{M}_{2g-2,\rho}$ provides an example of a ($g-1$)-dimensional family of structures 
in $\mathcal{M}_{2g-2,\rho}$ which are canonically branched but transverse to fibers of the Teichm\"uller map $\pi:\mathcal{M}_{2g-2,\rho}\to \mathcal T(S)$ (the critical direction should be tangent to the slice otherwise).
Contrast this with the case of a rational sphere $\Sigma_{(X,A)}$  (see Section \ref{sec_fromODEtoBPS} for details) 
which is a 1-dimensional family of 
structures in $\mathcal{M}_{2g-2,\rho}$ which are canonically branched and actually constitute a fiber of $\pi$. 
\end{remark}

\printbibliography

\end{document}